\RequirePackage{fix-cm}
\documentclass{amsart}      
\usepackage{graphicx}
\usepackage{amssymb,amsmath}
\usepackage{color}
\usepackage{mathabx}
\usepackage[all,cmtip]{xy}
\usepackage{tikz-cd}
\usepackage{stmaryrd}
\usepackage{hyperref}
\usepackage[]{algorithm2e}
\usepackage{algpseudocode}
\usepackage{pifont}
\usepackage[margin=1.5in]{geometry}

\newtheorem{thm}{Theorem}[section]
\newtheorem{prop}[thm]{Proposition}
\newtheorem{lem}[thm]{Lemma}
\newtheorem{cor}[thm]{Corollary}

\newcommand{\R}{\mathbb{R}}

\newcommand{\mathh}{\mathbb{H}}
\newcommand{\C}{\mathbb{C}}

\newcommand{\hatS}{\widehat{\mathcal{S}}}
\newcommand{\hatSop}{\widehat{\mathcal{S}}_{o}}
\newcommand{\mcSop}{\mathcal{S}_{o}}
\newcommand{\mcS}{\mathcal{S}}
\newcommand{\mcP}{\mathcal{P}}
\newcommand{\mcL}{\mathcal{L}}
\newcommand{\mrmH}{\mathrm{H}}

\newcommand{\St}{\mathrm{St}_2(\mathcal{V})}

\newcommand{\StOp}{\mathrm{St}_2^\ast(\mathcal{V})}
\newcommand{\StOpLoop}{\mathrm{St}^\ast_2(\mathcal{L}\C)}
\newcommand{\StOpAnti}{\mathrm{St}^\ast_2(\mathcal{A}\C)}
\newcommand{\Gr}{\mathrm{Gr}_2(\mathcal{V})}

\newcommand{\diffI}{\mathrm{Diff}^+(I)}

\begin{document}

\title{Shape Analysis of Framed Space Curves}


\author{Tom Needham}


\address{Department of Mathematics, The Ohio State University}


\maketitle

\begin{abstract}
In the elastic shape analysis approach to shape matching and object classification, plane curves are represented as points in an infinite-dimensional Riemannian manifold, wherein shape dissimilarity is measured by geodesic distance. A remarkable result of Younes, Michor, Shah and Mumford says that the space of closed planar shapes, endowed with a natural metric, is isometric to an infinite-dimensional Grassmann manifold via the so-called square root transform. This result facilitates efficient shape comparison by virtue of explicit descriptions of Grassmannian geodesics. In this paper, we extend this shape analysis framework to treat shapes of framed space curves. By considering framed curves, we are able to generalize the square root transform by using quaternionic arithmetic and properties of the Hopf fibration. Under our coordinate transformation, the space of closed framed curves corresponds to an infinite-dimensional complex Grassmannian. This allows us to describe geodesics in framed curve space explicitly. We are also able to produce explicit geodesics between closed, unframed space curves by studying the action of the loop group of the circle on the Grassmann manifold. Averages of collections of plane and space curves are computed via a novel algorithm utilizing flag means.

\end{abstract}

\section{Introduction}\label{sec:introduction}

The study of shape is of fundamental importance to problems in computer vision, object recognition, biomedical imaging and computer graphics. A common mathematical formalism for comparison of shapes of objects involves representing the objects as points in a metric space (the \emph{shape space}), where distance corresponds to some notion of shape dissimilarity. Classical approaches to shape analysis have represented objects by a finite number of  landmarks, whence the shape space is a finite-dimensional manifold \cite{DrydenBook,kendall1984shape}. There has recently been a large amount of work blending shape analysis with functional data analysis, wherein shapes are represented as equivalence classes of parameterized immersed manifolds \cite{bauer2014constructing,jermyn2012elastic,srivastava2011shape,younes2008metric}. The shape spaces in this framework are infinite-dimensional manifolds obtained as quotients of function spaces. Using this approach, a shape space can more generally be endowed with a Riemannian structure. This gives rise to a metric defined by geodesic distance as well as algorithms for statistical calculations such as Karcher means and Principal Component Analysis.

The purpose of this paper is to extend functional analysis-based shape comparison methods to the study of shapes of framed curves. A \emph{framed curve} is a curve in $\R^3$ endowed with a choice of adapted moving frame. There are two main motivations for doing so. First, framed curves have a wide variety of applications; for example, they are used in computer graphics and animation  \cite{bergou2010discrete,bertails2006super}, camera tracking \cite{goemans2004automatic}, modeling of protein folding \cite{hu2011discrete}, and they are central to the study of elasticity  \cite{kehrbaum1997elastic}. The second motivation is theoretical. By studying framed space curves, we are able to extend well-known work of Younes, Michor, Shah and Mumford on shapes of planar curves \cite{younes2008metric}. Roughly, their result says that the shape space of closed planar curves is isometric to an infinite-dimensional real Grassmann manifold. The simple geometry of the manifold leads to highly efficient shape comparisons, since geodesics in the Grassmannian are given explicitly. Our generalization says that the shape space of closed framed curves in $\R^3$ is isometric to a complex Grassmannian. Moreover, we show that the space of (unframed) closed curves in $\R^3$ is isometric to a quotient of the Grassmann manifold by the infinite-dimensional based loop group of $S^1$. The takeaway of our results is that calculations in these infinite-dimensional spaces of closed curves in $\R^3$ can be performed explicitly. A more precise description of the contributions of this paper follows below.

\subsection{Elastic Shape Analysis of Curves}\label{sec:elastic_shape_analysis}

In this section, we briefly outline the \emph{elastic shape analysis} approach to object matching. We will focus on shapes of curves in $\R^n$, although the general framework is quite flexible and has been extended to treat shapes of immersed surfaces \cite{jermyn2012elastic,kurtek2013landmark} and of curves in non-Euclidean manifolds \cite{LeBrigant2018,su2017square}. More details are available in several survey articles \cite{bauer2016use,srivastava2012advances} or the recent textbook \cite{srivastava2016functional}. 

For concreteness, we begin by considering shapes of closed curves in $\R^2$. Under the elastic shape analysis paradigm, one represents planar curves as points in an infinite-dimensional Riemannian manifold. The representation is achieved via a quotient construction; for practical purposes, we deal with parameterized curves, but we theoretically identify a pair of curves if they differ by a translation, a rotation, a scaling and/or a reparameterization. A point in the manifold of curve shapes is therefore the orbit of a parameterized curve under the action of the group of \emph{shape-preserving transformations}; i.e., the \emph{shape space of planar curves} is the quotient space
$$
\mathrm{Imm}(S^1,\R^2)/\left\{\begin{array}{c}\mathrm{transl.},\mathrm{rot.},\\
\mathrm{scal.},\mathrm{reparam.}\end{array}\right\},
$$
where $\mathrm{Imm}(S^1,\R^2)$ denotes the space of \emph{immersions} $\gamma: S^1 \rightarrow \R^2$ with $\gamma'(t) \neq 0$ for all $t$.

A Riemannian structure on the shape space is constructed by defining a Riemannian metric on the total space of immersed curves $\mathrm{Imm}(S^1,\R^2)$. If the metric is invariant under the  shape-preserving transformations, then it descends to a well-defined metric on the shape space.  There are many reasonable choices of metric on $\mathrm{Imm}(S^1,\R^2)$, and one typically considers Sobolev-type metrics of order at least one (due to the fact that geodesic distance vanishes with respect to the reparameterization-invariant $L^2$-metric \cite{michor2005vanishing}).  A popular class of metrics are those in the two parameter family $g^{a,b}$ of \emph{elastic metrics}  \cite{mio2007shape}, defined by
$$
g^{a,b}_\gamma(\mu,\mu) = \int_{S^1} a \left<D_s \mu, T\right>^2 + b \left< D_s \mu, N\right>^2 \; \mathrm{d}s,
$$
where $\gamma \in \mathrm{Imm}(S^1,\R^2)$, 
$$
\mu \in T_\gamma \mathrm{Imm}(S^1,\R^2) \approx C^\infty(S^1,\R^2),
$$
$T$ is the unit tangent to $\gamma$ and $N$ is its unit normal. Throughout the rest of the paper, we use 
$D_s = \frac{1}{\|\gamma'(t)\|} \frac{d}{dt}$ to denote derivative with respect to arclength and $\mathrm{d}s = \|\gamma'(t)\| \mathrm{d} t$ to denote arclength measure. 

A main inspiration for the present paper is the following well-known result of Younes, Michor, Shah and Mumford.
\begin{thm}[\cite{younes2008metric}]\label{thm:younes_thm}
The space 
$$
\mathrm{Imm}(S^1,\R^2)/\{\mathrm{transl.},\mathrm{rot.},\mathrm{scal.}\},
$$ 
endowed with the elastic metric $g^{a,a}$ is locally isometric to the Grassmann manifold $\mathrm{Gr}_2(C^\infty(S^1,\R))$ of two-dimensional planes in the vector space $C^\infty(S^1,\R)$, endowed with its canonical $L^2$ metric.
\end{thm}
The usefulness of this result lies in the fact that geodesics in the Grassmann manifold can be described explicitly by a straightforward extension of the finite-dimensional formula. This identification therefore facilitates a fast algorithm for shape comparisons. The proof of Theorem \ref{thm:younes_thm} relies completely on the identification of $\R^2$ with the complex plane, as the isometry is induced by the map $c(t) \mapsto \sqrt{c'(t)}$, with the square root taken pointwise, chosen so that the image is a continuous curve.

There have been several other papers building on this theme; that is, given an elastic metric on the space of plane curves, one constructs a coordinate transformation taking the elastic metric isometrically to a simplified metric \cite{bauer2014constructing,kurtek2018simplifying}. This procedure is sometimes referred to as \emph{flattening the metric}.
The most widely used of these is the \emph{square root velocity transform} (SRVT) of \cite{srivastava2011shape}, which flattens the elastic metric $g^{1/4,1}$. The SRVT generalizes to a transform for the space of immersed curves in $\R^n$, where immersion space is endowed with a straightforward generalization of the elastic metric $g^{1/4,1}$. 
 
A caveat is that, besides the square root transform of \cite{younes2008metric}, none of the metric-flattening transforms in the existing literature treat closed curves explicitly. As a result, most of the algorithms in the literature require, in each iterative step, an orthogonal  projection from the manifold of open curves to the submanifold of closed curves. Unfortunately, it is unclear how the square root transform and its resulting explicit geodesics for closed curves would directly generalize to curves in higher dimension, due to its reliance on complex arithmetic. We will show that if one considers \emph{framed} curves in $\R^3$, then the pointwise complex arithmetic can be replaced with quaternionic arithmetic, and there is an appropriate generalization of Theorem \ref{thm:younes_thm} which treats closed curves explicitly.

\subsection{Main Contributions and Outline of the Paper}

We begin Section \ref{sec:shape_analysis_of_open_framed_curves} with definitions of our spaces of framed curves. A family of elastic metrics is defined on the space of framed curves and we show that the family is a natural generalization of the elastic metrics for plane curves (a different generalization than the one used in the literature on SRVT for curves in $\R^n$). We then utilize quaternionic algebra to give a coordinate transformation which flattens a particular choice of framed curve elastic metric. This construction was described in  previous articles \cite{needham2017knot,needham2017kahler}, where the focus was on symplectic geometry and theoretical applications. The exposition provided here is focused on algorithms and applications to shape analysis. We show that the transformation takes our elastic metric isometrically to an $L^2$ metric on the target space (Theorem \ref{thm:pullback_metric}). Using this transformation, we give an explicit formula for geodesics in the space of open framed curves.

Section \ref{sec:shape_analysis_of_closed_framed_curves} specializes to closed framed curves. We show how the coordinate transformation described above restricts to identify the shape space of closed framed curves with a complex Grassmannian (Theorem \ref{thm:grassmannian}). An explicit formula for geodesics in the complex Grassmannian is provided.

In Section \ref{sec:frame_twisting}, we introduce another group action which is genuinely unique to the framed curve setting. This is the action of $C^\infty(S^1,S^1)$ on framed curve space by \emph{frame twisting}. The space of \emph{unframed} curves in $\R^3$ can be viewed as the quotient of framed curve space by this group action. Theorem \ref{thm:frame_twisting} shows that geodesics in framed curve space which are horizontal with respect to this action can be determined explicitly. As a consequence, we are able to find explicit geodesics for open and closed (unframed) curves in $\R^3$ with respect to our metric.

The paper concludes in Section \ref{sec:applications} with several examples of geodesics between framed and unframed curves. We also apply our theory to shape analysis of supercoiled circular DNA molecules, with shape data from the experiment described in \cite{irobalieva2015structural}. In this section, we also introduce a new algorithm for averaging plane curves and space curves, based on the notion of a flag mean for points sampled from a Grassmannian \cite{draper2014flag}.

\section{Shape Analysis of Open Framed Curves}\label{sec:shape_analysis_of_open_framed_curves}

In this section we describe the structure of the shape space of open framed curves and our algorithm for producing geodesics between framed curves.

\subsection{Spaces of Framed Paths}

\subsubsection{The Preshape Space}

A \emph{framed space curve} is a pair $(\gamma,V)$ of smooth maps $\gamma,V: I \rightarrow \R^3$ from a closed interval $I \subset \R$ such that $\gamma$ is an immersion and $V$ is a unit normal vector field to $\gamma$. We fix the choice of domain $I=[0,2]$ as a convenient normalization; the utility of this particular choice will become clear when we pass to the subspace of closed curves. Our end goal is to study \emph{shapes} of framed curves, which are equivalence classes with respect to the shape-preserving group actions of translation, scaling, rotation and reparameterization. 

The simplest invariance to treat is the action of $\R^3$ by translation, since the quotient by $\R^3$ is concretely realized by representing all framed curves $(\gamma,V)$ with a fixed basepoint $\gamma(0)=\vec{0}$. The resulting quotient space $\hatS_{o}$ will be referred to as the \emph{preshape space of (open) framed paths}.

\subsubsection{Manifold Structure}

In this paper, infinite-dimensional manifolds will belong to the Nash-Moser category of tame Fr\'{e}chet manifolds \cite{hamilton1982inverse}. For a finite-dimensional manifold $M$, we use the simplified notation $\mathcal{P}M = C^\infty(I,M)$ and $\mathcal{L}M = C^\infty(S^1,M)$ for the \emph{path space} and \emph{loop space} of $M$, respectively. 

It is shown in \cite{needham2017kahler} that $\hatS_{o}$ is an infinite-dimensional tame Fr\'{e}chet manifold by identifying it with $\mathcal{P}(\mathrm{SO}(3) \times \R^+)$ using the map
\begin{equation}\label{eqn:map_to_SO3}
(\gamma,V) \mapsto \left(\left(T,V,T\times V\right),\|\gamma'\|\right),
\end{equation}
where $\|\cdot\|$ will always denote the Euclidean norm and $T=\gamma'/\|\gamma'\|$ denotes the \emph{unit tangent curve} of $\gamma$. By fixing the standard basis, we identify $\mathrm{SO}(3)$ with the space of $3 \times 3$ matrices with orthonormal columns.  Each entry $T$, $V$ and $T \times V$ in the image of \eqref{eqn:map_to_SO3} is then considered as a path of column vectors. The map \eqref{eqn:map_to_SO3} is a bijection: the only information lost by the map is the basepoint $\gamma(0)$, and this is accounted in the definition of $\hatSop$.

\subsubsection{The Shape Space}

The finite-dimensional  Lie groups $\R^+$ and $\mathrm{SO}(3)$ act on $\hatSop$ by scaling (of the base curve) and rotation, respectively, while the infinite-dimensional Lie group $\diffI$ of orientation-preserving diffeomorphisms of $I$ acts by precomposition: for $\rho \in \diffI$ and $(\gamma,V) \in \hatSop$, the action is given by the formula
$$
\rho \cdot (\gamma,V) = (\gamma \circ \rho, V \circ \rho).
$$
We refer to this as the \emph{reparameterization action}. The \emph{shape space of open framed curves} is then defined to be the quotient space
$$
\mcSop=\hatSop/(\R^+ \times \mathrm{SO}(3) \times \diffI).
$$

\subsection{Metrics on the Shape Space}

\subsubsection{Framed Path Elastic Metrics}

In order to perform statistical analysis on collections of framed curves, we wish to introduce a Riemannian metric on the shape space $\mcSop$. Following the elastic shape analysis framework outlined in the introduction, this is accomplished by first defining a metric on the preshape space $\hatSop$ which is invariant under scaling, rotation and reparameterization, so that it descends to a well-defined metric on the quotient shape space. 

It is a straightforward exercise to verify that the tangent space to a base point $(\gamma,V) \in \hatSop$ consists of smooth variations $(\nu,W)$ satisfying the constraints
\begin{align*}
\nu(0)&=\vec{0} \\
\left< \nu' , V \right> + \left<\gamma', W\right> = \left< W,V\right> &= 0.
\end{align*} 
Such a variation decomposes at each point into four components: bending toward $V$, bending toward $T \times V$, stretching of $\gamma$, and rotation of $V$ around $T$. Inspired by the plane curve elastic metrics introduced in Section \ref{sec:elastic_shape_analysis}, we define the \emph{framed curve elastic metric} with parameters $a,b,c,d$ by
\begin{align*}
g^{a,b,c,d}_{(\gamma,V)} \left((\nu,W),(\nu,W)\right)  &= \int_I  a \left<D_s \nu, V\right>^2  + b \left<D_s \nu,T \times V\right>^2 \\
& \hspace{.4in} + c \left<D_s \nu,T\right>^2  + d \left<W ,T \times V\right>^2   \; \mathrm{d}s.
\end{align*}
The parameters control the weights of the four types of deformations described above. Since the derivatives and the measure are with respect to arclength, these metrics are invariant under $\diffI$. Moreover, they are translation and rotation-invariant.  Finally, they scale with homotheties of the base curve and can therefore be made scale invariant by dividing by the total length of $\gamma$; we instead opt to treat this invariance by working primarily in the submanifold of curves of fixed length.

\subsubsection{The Submanifold of Planar Curves}

The space $C^\infty(I,\R^2)/\R^2$ of based planar curves embeds naturally into $\hatSop$ as the submanifold of framed curves $(\gamma,V)$ such that the image of $\gamma$ lies in the $xy$-plane and $V$ is the oriented unit normal vector to $\gamma$. The tangent variations to this submanifold have no bending component in the $T \times V$-direction and no twisting component. These terms therefore vanish in the restriction of $g^{a,b,c,d}$ to this submanifold and it follows that the induced metric on the submanifold is the planar elastic metric $g^{a,c}$. The framed curve elastic metrics therefore give a natural generalization of the plane curve elastic metrics which is an alternative to the standard generalization to curves in $\R^n$ \cite{srivastava2011shape}.

 \subsubsection{Fixing a Parameter Choice}
 
In this paper, we focus on the particular parameter choice $a=b=c=d=1$. Let $g^\mcS$ be the metric given by the formula
$$
g^\mathcal{S}_{(\gamma,V)} \left((\nu_1,W_1),(\nu_2,W_2))\right) = \frac{1}{4} \int_I  \left<D_s \nu_1, D_s \nu_2\right>  + \left<W_1,T \times V\right>\left<W_2,T \times V\right>  \; \mathrm{d}s.
$$
This metric will be used for all of our shape spaces, hence the generic superscript in our notation. A simple calculation shows that 
$$
g^\mcS = \frac{1}{4} g^{1,1,1,1}.
$$
We rescale by a factor of $4$ as a matter of convenience; this scaling will disappear when we introduce a change of coordinates in the following subsection.

When restricting to the submanifold of planar curves, the induced metric is a constant multiple of the planar elastic metric $g^{1,1}$, and the work in this paper can therefore be seen as a direct extension of the results of \cite{younes2008metric}.

\subsection{Quaternionic Coordinates for Framed Paths}

\subsubsection{The Frame-Hopf Map}

Let $\mathbb{H}=\mathrm{span}_\R \{1,\textbf{i}, \textbf{j}, \textbf{k}\}$ denote the skew-field of quaternions. A quaternion
$$
q=q_0 + q_1 \textbf{i} + q_2 \textbf{j} + q_3 \textbf{k} \in \mathbb{H}
$$
has conjugate denoted
$$
\overline{q} = q_0 - q_1 \textbf{i} - q_2 \textbf{j} - q_3 \textbf{k}.
$$
We will abuse notation slightly and also denote quaternion-valued paths and loops by $q=q(t)$.

It is well known that a framed curve $(\gamma,V)$ can be represented as a path in the quaternions $q \in \mathcal{P}\mathh$ which is unique up to a global choice of sign (e.g., \cite{dichmann1996hamiltonian,hanson2005visualizing}). More precisely, there is a smooth double covering map $\mathrm{H}:\mathcal{P}\mathh^\ast \rightarrow \hatSop$ defined explicitly by
\begin{equation}\label{eqn:hopf_map_formula}
\mathrm{H}(q) =(\gamma,V) = \left(\int \overline{q}\textbf{i} q \; \mathrm{d}t, \frac{\overline{q}\textbf{j}q}{\|q\|_\mathh^2}\right) 
\end{equation}
with $\mathrm{H}(q_1)=\mathrm{H}(q_2)$ if and only if $q_1=\pm q_2$, which we refer to as the \emph{frame-Hopf map}. In \eqref{eqn:hopf_map_formula}, all quaternionic arithmetic in the formula is understood to be performed pointwise on the curves. It is easy to check that the quaternionic paths $\overline{q}\textbf{i}q$ and $\overline{q}\textbf{j}q$ are purely imaginary, so that they can be naturally identified with paths in $\R^3$. The norm $\|\cdot\|_\mathh$ is the Euclidean norm on $\mathh \approx \R^4$, and is also applied pointwise. The integral symbol in the first coordinate denotes the antiderivative based at $\vec{0}$.  It is also straightforward to check that $V$ defines a normal vector field to $\gamma$. Putting all of this together, we conclude that the frame-Hopf map is well-defined.

\subsubsection{The Classical Hopf Map}

We call $\mrmH$ the frame-Hopf map in reference to its relationship to the well known anti-homomorphic double-covering 
\begin{equation}\label{eqn:homomorphism}
h:\mathrm{SU}(2) \rightarrow \mathrm{SO}(3). 
\end{equation}
The double-covering $h$ can be realized by identifying $\mathrm{SU}(2)$ with the $3$-sphere of unit quaternions $S^3 \subset \mathh$ via
$$
\left(\begin{array}{cc}
u & v \\
-\overline{v} & \overline{u} \end{array}\right) \leftrightarrow u+v \textbf{j},
$$
where $u$ and $v$ are complex numbers satisfying $|u|^2+|v|^2=1$. Then $h:S^3 \rightarrow \mathrm{SO}(3)$ is the map
$$
q \mapsto (\overline{q}\textbf{i}q,\overline{q}\textbf{j}q,\overline{q}\textbf{k}q).
$$
Each entry in the image is a purely imaginary quaternion,  and can therefore be identified with a column vector in $\R^3$. Then each column has unit norm and the columns are pairwise orthogonal, so the image lies in $\mathrm{SO}(3)$. Up to multiplication by a constant, this map is an isometry with respect to the natural invariant metrics on $\mathrm{SU}(2)$ and $\mathrm{SO}(3)$. See, e.g., \cite{gelfand1963representations} for more details.

\subsubsection{Local Isometry Theorem}

Let $\left<\cdot,\cdot\right>_\mathh$ denote the Euclidean inner product on $\mathh \approx \R^4$. In quaternionic coordinates, this is given by 
$$
\left<q_0,q_1\right>_\mathh = \mathrm{Re} \, q_0 \overline{q_1}.
$$ 
We denote the standard $L^2$ metric on the vector space $\mathcal{P}\mathh$ by
$$
\left<q_1,q_2\right>_{L^2} = \int_I \left<q_1,q_2\right>_\mathh \; \mathrm{d}t
$$
and its induced norm by $\|\cdot \|_{L^2}$.

The following theorem relates the complicated metric $g^\mcS$ on framed path space to this simple $L^2$-metric. This result is stated  in \cite{needham2017kahler}, but not proved. We include a proof in the appendix (Section \ref{sec:proof_pullback}).

\begin{thm}\label{thm:pullback_metric}
The pullback of $g^\mcS$ by $\mrmH$ satisfies
$$
\mrmH^\ast g^\mcS = g^{L^2}.
$$
\end{thm}

This motivates us to define a shape similarity metric on $\hatSop$ by assigning a pair of framed curves $(\gamma_j,V_j)$, $j=0,1$, the distance
\begin{equation}\label{eqn:distance_on_preshapes}
\min \{\|q_0 - q_1\|_{L^2} \},
\end{equation}
where the minimum is taken over the four possible combinations of lifts of of the $(\gamma_j,V_j)$. Theorem \ref{thm:pullback_metric} implies that this is equal to the geodesic distance between the framed curves with respect to the metric $g^\mathcal{S}$. 

\subsection{Modding out Scaling and Rotation}

\subsubsection{Scale-Invariance}

To explicitly treat scale-invariance for framed curve shapes, we can preprocess and consider framed curves $(\gamma,V)$ such that the base curve $\gamma$ has fixed length $2$ (once again, a convenient normalization). This is the space
$$
\hatSop/\R^+ \approx \left\{(\gamma,V) \in \hatSop \mid \int_I \|\gamma'\| \; \mathrm{d}t = 2\right\},
$$
which is a codimension-1 submanifold of $\hatSop$ (this follows by an easy application of Hamilton's implicit function theorem \cite[Section III, Theorem 2.3.1]{hamilton1982inverse}). The frame-Hopf map restricts to give an isometric double covering 
$$
\mathrm{H}:S^\ast_{\sqrt{2}} \rightarrow \hatSop/\R^+,
$$ 
where $S_{\sqrt{2}}$ is the radius-$\sqrt{2}$ \emph{Hilbert sphere}
$$
S_{\sqrt{2}} = \left\{ q \in \mathcal{P}\mathh \mid \|q\|_{L^2}^2 = 2 \right\}
$$
and
$$
S_{\sqrt{2}}^\ast = S_{\sqrt{2}} \cap \mathcal{P}\mathh^\ast.
$$
Geodesic distance between points $(\gamma_0,V_0)$ and $(\gamma_1,V_1)$ in the submanifold $\hatSop/\R^+$ is therefore given by 
$$
\sqrt{2} \cdot \min\{\arccos \left<q_0,q_1\right>_{L^2}\},
$$
where the minimum is taken over lifts $q_j$ of the $(\gamma_j,V_j)$ and $\sqrt{2} \cdot \arccos \left<q_0,q_1\right>_{L^2}$ is great-circle distance in the sphere $S_{\sqrt{2}}$.

\subsubsection{Rotation-Invariance}

Rotation-invariance for framed curves is naturally built into this framework as well. The group $\mathrm{SU}(2)\approx S^3$ acts on $\mcP \mathh^\ast$ by pointwise multiplication. The frame-Hopf map has the following equivariance property for all $A \in \mathrm{SU}(2)$ and $q \in \mcP \mathh^\ast$:
$$ 
\mathrm{H}(q \cdot A) = h(A) \cdot \mathrm{H}(q).
$$
The action on the righthand side is the $\mathrm{SO}(3)$-action on $\hatSop$ by pointwise rotation, where the matrix $h(A)$ is the image of $A$ under the map \eqref{eqn:homomorphism}. This equivariance descends to the restricted map $\mathrm{H}:S_{\sqrt{2}} \rightarrow \hatSop/\R^+$.

Recall that $g^\mcS$ is rotation-invariant. This implies that it descends to a well-defined metric on the quotient $\hatSop/(\R^+ \times \mathrm{SO}(3))$. By the above discussion, geodesic distance between the equivalence classes of framed curves $(\gamma_0,V_0)$ and $(\gamma_1,V_1)$ in this quotient is given by
\begin{equation}\label{eqn:modding_out_rotation}
\sqrt{2} \cdot \min\{\arccos \left< q_0 \cdot A_0,  q_1 \cdot A_1\right>_{L^2}\},
\end{equation}
where the minimum is taken over quaternionic lifts of the $(\gamma_j,V_j)$ and over elements $A_j$ of the compact Lie group $\mathrm{SU}(2)$. Optimal alignment over rotations is treated explicitly by the following proposition. The proof is in Section \ref{sec:proof_rotations}.

\begin{prop}\label{prop:rotation_optimization}
Let $q_0, q_1 \in S_{\sqrt{2}}$ be quaternionic paths which are not $L^2$-orthogonal. The minimum distance between the $\mathrm{SU}(2)$-orbits of $q_0$ and $q_1$ is realized by
$$
\sqrt{2} \cdot \arccos \left<q_0, q_1\cdot \widehat{A}\right>_{L^2},
$$
where $\widehat{A}$ is the normalization of the quaternion
$$
\int_I \overline{q_0} q_1 \; \mathrm{d}t.
$$
\end{prop}

\subsection{Reparameterizations and Geodesics in Shape Space}

\subsubsection{The $\diffI$-action in Quaternionic Coordinates}

A diffeomorphism  $\rho \in \diffI$ acts on $q \in \mcP \mathh^\ast$ via
$$
\rho \cdot q = \sqrt{\rho'}(q \circ \rho).
$$
It is straightforward to check that the map $\mathrm{H}$ is equivariant with respect to this action on $\mcP \mathh^\ast$ and the reparameterization action on $\hatSop$; i.e.,
$$
\mathrm{H}\left(\sqrt{\rho'}(q \circ \rho)\right) = \mathrm{H}(q)\circ \rho.
$$

Geodesic distance between equivalence classes of framed curves $(\gamma_1,V_1)$ and $(\gamma_2,V_2)$ in $\hatSop/\diffI$ is therefore given by
\begin{equation}\label{eqn:distance_diff_1}
\inf \{\|\rho_0 \cdot q_0 - \rho_1 \cdot q_1\|_{L^2}\},
\end{equation}
where the infimum is taken over lifts $q_j$ of the $(\gamma_j,V_j)$ and over all $\rho_j \in \diffI$. The fact that $\diffI$ acts by $L^2$ isometries implies that this expression can be simplified to
\begin{equation}\label{eqn:distance_diff}
\inf\{\|q_0 - \rho \cdot q_1\|_{L^2}\},
\end{equation}
where the infimum is taken over lifts and $\rho \in \diffI$. Moreover, the $\diffI$-action preserves the submanifold of fixed length curves (or the Hilbert sphere, in quaternionic coordinates) and commutes with the rotation action of $\mathrm{SO}(3)$. Geodesic distance in the shape space $\mathcal{S}_o$ is then given by
\begin{equation}\label{eqn:distance_shape_space}
\inf \{\sqrt{2} \cdot \arccos \left<q_0, \rho \cdot q_1 \cdot A\right>_{L^2}\},
\end{equation}
with the infimum taken over all quaternionic lifts $q_j$, all $\rho \in \diffI$ and all $A \in \mathrm{SU}(2)$.

\subsubsection{Existence of Optimal Reparameterizations}\label{sec:existence_of_reparam}

It is natural to ask whether the infima in the distance formulas \eqref{eqn:distance_diff} and \eqref{eqn:distance_shape_space} are realized by smooth reparameterizations $\rho$. In fact, they are not; an example of Younes et.\ al.\ appearing in \cite[Section 4.2]{younes2008metric} in the context of planar curves also applies here, since our metric restricts to the Younes metric on the submanifold of plane curves. 

The corresponding question in the SRVT setting has been rigorously studied in several recent articles and is known to have a positive answer in lower regularity classes. It is shown in \cite{bruveris2016optimal} that the optimal matching between $C^1$ curves is realized by an absolutely continuous reparameterization. The case of piecewise linear curves is considered in \cite{lahiri2015precise}, where it is shown that the optimal reparameterization is realized by a PL curve. This result was extended to more general elastic metrics on PL plane curves in \cite{kurtek2018simplifying}. 

Precise results in lower regularity classes for the framed curve setting will be the subject of future work. In this article, we take the practical point of view that optimal reparameterizations will be approximated numerically. The step of optimizing over $\diffI$, requires an approximation of 
\begin{equation}\label{eqn:approximate_over_diff}
\inf_\rho  \int \left\| q_0 - \sqrt{\rho'} q_1 \circ \rho \right\|_\mathrm{H}^2 \; \mathrm{d}t,
\end{equation}
for $q_j \in \mathcal{P}\mathrm{H}^\ast \approx \mathcal{P}\left(\R^4 \setminus \{0\}\right)$. Fortuitously, the identification of $\|\cdot\|_\mathbb{H}$ with the Euclidean norm in $\R^4$ implies that this is exactly the optimization problem which appears in the SRVT setting, meaning that existing, highly efficient algorithms can be directly applied. We will use the dynamic programming algorithm of \cite{mio2007shape} to approximate solutions of \eqref{eqn:approximate_over_diff}.

\subsubsection{Geodesics Between Open Framed Curves}

Geodesics in the shape space $\mathcal{S}_o$ are computed as follows. Given framed curves $(\gamma_0,V_0)$ and $(\gamma_1,V_1)$, we preprocess as necessary to ensure $\mathrm{length}(\gamma_j)=2$. The quaternionic lifts of the framed curves $q_0$ and $q_1$ are then elements of the Hilbert sphere $S_{\sqrt{2}}$. Numerical algorithms for lifting a framed curve to a quaternionic representation can be found in \cite{hanson2005visualizing}. We choose our lifts to minimize $L^2$ distance. We then compute an optimal rotation via Proposition \ref{prop:rotation_optimization} and an approximate optimal reparameterization via the dynamic programming algorithm of \cite{mio2007shape}. This procedure is iterated until a stopping condition is met (e.g., decreases in geodesic distance fall below a fixed threshold, or a predetermined number of iterations is reached) to produce a final reparameterization $\widehat{\rho}$ and rotation $\widehat{A}$. Let $\widehat{q}_1 = \widehat{\rho} \cdot q_1 \cdot \widehat{A}$. The geodesic joining the equivalence classes of the framed curves $(\gamma_j,V_j)$ in the shape space is realized by applying the frame-Hopf map $\mathrm{H}$ pointwise to the spherical interpolation
\begin{equation}\label{eqn:spherical_geodesic}
q_u = \frac{\sin((1-u)\theta)}{\sin \theta} q_0 + \frac{\sin(u \theta)}{\sin \theta} \widehat{q}_1,
\end{equation}
with $u \in [0,1]$ and with $\theta$ denoting geodesic distance between $q_0$ and $\widehat{q}_1$. Several examples of geodesics are provided in Section \ref{sec:applications}.

We note that $S_{\sqrt{2}}^\ast$ is not geodesically complete; there are $q_0$ and $\widehat{q}_1$ for which the geodesic given by \eqref{eqn:spherical_geodesic} passes through $S_{\sqrt{2}} \setminus S_{\sqrt{2}}^\ast$. This does not cause an issue from a practical point of view, as geodesic distance is still a reasonable measure of shape dissimilarity between the shapes and the frame-Hopf map $\mathrm{H}$ is still well-defined at those points---they correspond geometrically to singular curves along the geodesic. Moreover,  it holds generically that geodesics stay in $S_{\sqrt{2}}^\ast$, in the sense that a generic homotopy of curves in $\R^4$ will not pass through the origin.

\section{Shape Analysis of Closed Framed Curves}\label{sec:shape_analysis_of_closed_framed_curves}

We now turn to a similar description of the geometry of the space of closed framed curves. Unlike the SRVT setting, geodesics in the submanifold of closed curves can be described explicitly. This is facilitated by a natural identification of closed curve space with an infinite-dimensional Grassmannian, generalizing Theorem \ref{thm:younes_thm}. The construction of this identification was described in \cite{needham2017kahler}, but is summarized here for the convenience of the reader.

\subsection{The Space of Closed Framed Curves}

\subsubsection{The Pre-Shape Space of Framed Loops}

A \emph{framed loop} is a pair $(\gamma,V)$ of smooth maps from the circle $S^1$ into $\R^3$ such that $\gamma$ is an immersion and $V$ is a (periodic) unit normal vector field to the image of $\gamma$. In order to identify the collection of framed loops with a subset of the shape space of framed paths, we identify $S^1$ with the quotient $[0,2]/0\sim 2$. It will still be convenient to mod out by translations and we denote by $\hatS_c$ the \emph{preshape space of} \emph{(closed) framed loops}, which is the set of framed loops satisfying $\gamma(0)=\vec{0}$. Just as in the framed paths case, we will treat remaining shape-preserving group actions separately.

\subsubsection{Manifold Structure}\label{sec:closed_manifold_structure}

The map \eqref{eqn:map_to_SO3} restricts to $\hatS_c$ to give an injective map into $\mcL (\mathrm{SO}(3) \times \R^+)$, but it is no longer surjective; for example, any constant loop is not in its image. The image is a codimension-3 submanifold of $\mcL(\mathrm{SO}(3) \times \R^+)$ and this can be used to endow $\hatS_c$ with a manifold structure \cite{needham2017kahler}. It can be deduced from the fact that the fundamental group of $\mathrm{SO}(3)$ is $\mathbb{Z}/2\mathbb{Z}$ that the manifold $\hatS_c$ has two connected components. The component that a framed $(\gamma,V)$ curve belongs to is determined by the topological linking number of $\gamma$ with a small pushoff along the $V$--direction, measured modulo-2 (a more detailed explanation is given in \cite{needham2017kahler,needham2016grassmannian}).

\subsubsection{The Shape Space of Framed Loops}

Similar to the case of open framed curves, the Lie groups $\R^+$, $\mathrm{SO}(3)$ and $\mathrm{Diff}^+(S^1)$ act on $\hatS_c$ by scaling, rotation and reparameterization, respectively. There is a remaining group action by which it will be convenient to quotient. This is the action of $S^1$ by \emph{global frame twists} (as opposed to local frame twisting, which is treated in section \ref{sec:frame_twisting}). An element $\theta \in S^1$ acts on a framed loop $(\gamma,V)$ by rotating each vector $V(t)$ by the angle $\theta$ in the plane normal to $\gamma(t)$ with respect to the right-hand rule. We define the \emph{shape space of framed loops} to be the quotient space
$$
\mcS_c=\hatS_c/(\R^+ \times \mathrm{SO}(3) \times S^1 \times \mathrm{Diff}^+(S^1)).
$$

\subsection{Complex Coordinates for Framed Loops}

\subsubsection{Antiloop Space}

We wish to determine the subset of $\mathcal{P}\mathh^\ast$ which covers $\hatS_c$ via the frame-Hopf map. We first note that closure is not a necessary condition for a quaternionic curve to be mapped to a smoothly closed framed loop. For example, the open quaternionic curve
\begin{align*}
q(t)&=\cos(\pi t/2) + \textbf{i} \sin (\pi t/2) + \\
&\hspace{0.6in} \textbf{j} \cos (\pi t/2) - \textbf{k} \sin (\pi t/2)
\end{align*}
maps to a closed framed loop under the frame-Hopf map. We introduce the \emph{anti-loop space} 
$$
\mathcal{A}\mathh^\ast = \{q \in \mathcal{P}\mathh^\ast \mid q^{(k)}(2) = - q^{(k)}(0)\, \forall k\}.
$$
A necessary condition for $q \in \mathcal{P}\mathh^\ast$ to correspond to a smoothly closed framed loop is that 
\begin{equation}\label{eqn:necessary_condition}
q \in \mcL \mathh^\ast  \sqcup \mathcal{A}\mathh^\ast. 
\end{equation}
The necessity of this condition follows easily from the double-covering property of the frame-Hopf map. The disjoint union here corresponds to the fact that $\hatS_c$ has two path components.

\subsubsection{Closure Condition}

Condition \eqref{eqn:necessary_condition} on $q$ is not sufficient for $\mathrm{H}(q)$ to be a framed loop, due to the codimensionality of the closure condition discussed in Section \ref{sec:closed_manifold_structure}. To describe the closure condition in quaternionic coordinates, it is useful to identify $\mathh$ with $\C^2$ via
$$
q_0 + q_1 \textbf{i} + q_2 \textbf{j} + q_3 \textbf{k} \leftrightarrow (q_0 + q_1 i, q_2 + q_3 i).
$$
This gives an identification of inner product spaces $(\mathh,\left<\cdot,\cdot\right>_\mathh)$ and $(\C^2,\mathrm{Re}\left<\cdot,\cdot\right>_{\C^2})$, where $\left<\cdot,\cdot\right>_{\C^2}$ is the standard Hermitian inner product on $\C^2$ and this extends to an identification of inner product spaces $\mathcal{P}\mathh$ and $\mathcal{P}\C^2$ with their induced $L^2$ inner products. We will convert between quaternionic and complex notation freely throughout the rest of the paper.

It will be useful to describe the frame-Hopf map in complex coordinates. For $(z,w)$ in $\mathcal{P}\C^2$, $\mathrm{H}(z,w)=(\gamma,V)$ is given by the formulas
\begin{align}
\gamma &= \int \left(|z|^2-|w|^2,2\mathrm{Im}(z\overline{w}),2\mathrm{Re}(z\overline{w})\right)\;\mathrm{d}t, \label{eqn:hopf_map_complex_coords} \\
V &= \frac{1}{(|z|^2 + |w|^2)}\left(2\mathrm{Im}(zw), \right.  \nonumber \\
&\hspace{.4in} \left.\mathrm{Re}(z^2+w^2),\mathrm{Im}(-z^2+w^2)\right), \nonumber
\end{align}
where the integral in the formula for $\gamma$ denotes the antiderivative based at zero. Using these explicit formulas, the following lemma is immediate.

\begin{lem}\label{lem:fundamental_lemma}
Let $q=(z,w) \in \mathcal{P}\mathh^\ast$ with $\mathrm{H}(q)=(\gamma,V)$.
\begin{itemize}
\item[(a)] $(\gamma,V)$ is a smoothly closed framed curve if and only if $q$ is smoothly closed or anticlosed and $z$ and $w$ are $L^2$-equinorm and orthogonal.
\item[(b)] The length of $\gamma$ is given by $\|q\|^2_{L^2}=\|z\|_{L^2}^2+\|w\|_{L^2}^2$. 
\end{itemize}
\end{lem}

\subsubsection{Stiefel Manifolds}

Throughout the rest of the paper, we generically use $\mathcal{V}$ for $\mathcal{L}\C$ or $\mathcal{A}\C$ (defined analogously to $\mathcal{A}\mathh^\ast$). The \emph{Stiefel manifold of orthonormal $2$-frames in $\mathcal{V}$} is the space
\begin{align*}
\St &= \{(z,w) \in \mathcal{V}^2 \mid \\
&\hspace{.2in} \|z\|_{L^2} = \|w\|_{L^2} = 1, \, \left<z,w\right>_{L^2}= 0\}.
\end{align*}
By a slight abuse of notation, we use $\left<\cdot,\cdot\right>_{L^2}$ for the standard $L^2$ inner product on $\mathcal{V}$. We are particularly interested in the open submanifold $\StOp$ consisting of $(z,w)$ in $\St$ such that $(z(t),w(t)) \neq \vec{0}$ for all $t$. Lemma \ref{lem:fundamental_lemma} immediately implies the following corollary.

\begin{cor}\label{cor:stiefel_manifold}
The frame-Hopf map restricts to an isometric double covering
$$
\StOpLoop \sqcup \StOpAnti \rightarrow \hatS_c/\R^+
$$
with respect to the $L^2$ metric and $g^\mathcal{S}$.
\end{cor}

\subsubsection{Grassmann Manifolds}

As we have seen, the frame-Hopf map is equivariant with respect to the rotation actions of $\mathrm{SU}(2)$ and $\mathrm{SO}(3)$. Moreover, the action of the diagonal circle $\mathrm{U}(1) \subset \mathrm{U}(2)$ on $\mathcal{P}\mathh^\ast$ corresponds to the frame twisting action of $S^1$ on framed path space in the sense that 
$$
\mathrm{H}\left(q \cdot \left(\begin{array}{cc}
\theta & 0 \\
0 & \theta \end{array}\right)\right) = (2\theta) \cdot \mathrm{H}(q).
$$
This motivates us to quotient by the full unitary group to obtain the \emph{Grassmann manifolds}
$$
\mathrm{Gr}_2(\mathcal{V})) = \mathrm{St}_2(\mathcal{V})/\mathrm{U}(2)
$$
and
$$
\mathrm{Gr}^\ast_2(\mathcal{V}) = \mathrm{St}^\ast_2(\mathcal{V})/\mathrm{U}(2).
$$

Corollary \ref{cor:stiefel_manifold} and the equivariance property described above combine to give the following theorem.

\begin{thm}[\cite{needham2017kahler}]\label{thm:grassmannian}
The frame-Hopf map induces an isometry
$$
\mathrm{Gr}^\ast_2(\mathcal{L}\C) \sqcup \mathrm{Gr}^\ast_2(\mathcal{A}\C)\rightarrow \hatS_c/(\R^+ \times \mathrm{SO}(3) \times S^1)
$$
with respect to the induced $L^2$ metric and $g^\mathcal{S}$.
\end{thm}

\subsection{Geodesics for Closed Framed Curves}

\subsubsection{Grassmannian Geodesics}\label{sec:grassmannian_geodesics}

Theorem \ref{thm:grassmannian} implies that geodesics in the quotient 
$$
\hatS_c/(\R^+ \times \mathrm{SO}(3) \times S^1)
$$
can be realized as  images of geodesics in the Grassmannian under the frame-Hopf map. It is therefore a major benefit of this approach that geodesics can be described explicitly in these manifolds via the method of \emph{Neretin geodesics}, adapted to the complex setting \cite{neretin2001jordan}. 

We represent an element of $\Gr$ as the $\mathrm{U}(2)$-orbit of a point $(z,w) \in \St$. We denote the orbit of $(z,w)$ by $[z,w]$; note that this can be interpreted geometrically as the $2$-plane in $\mathcal{V}$ spanned by $\{z,w\}$. Given two points $(z_0,w_0)$ and $(z_1,w_1)$ of the Steifel manifold $\St$, we construct a geodesic $(z_u,w_u)$, $u \in [0,1]$ joining their orbits in $\Gr$ as follows:
\begin{enumerate}
\item Compute the singular value decomposition of the projection map $[z_0,w_0] \rightarrow [z_1,w_1]$ between $2$-planes. From the SVD we obtain  a new orthonormal basis $\left(\widetilde{z}_j,\widetilde{w}_j\right)$ for $[z_j,w_j]$ such that the projection map is given by $\widetilde{z}_0 \mapsto \lambda_z \widetilde{z}_1$ and $\widetilde{w}_0 \mapsto \lambda_w \widetilde{w}_1$, where $0 \leq \lambda_z,\lambda_w \leq 1$ are singular values. The choice of basis corresponds to registration over $\mathrm{SO}(3) \times S^1$ in  framed loop space.
\item The Jordan angles of the $2$-planes are $\theta_z = \arccos(\lambda_z)$ and $\theta_w = \arccos(\lambda_w)$. 
\item The geodesic in $\Gr$ is realized by the interpolations
\begin{align*}
z_u &= \frac{\sin ((1-u) \cdot \theta_z) \widetilde{z}_0 + \sin(u \cdot \theta_z) \widetilde{z}_1}{\sin \theta_z} \\
w_u & = \frac{\sin ((1-u) \cdot \theta_w) \widetilde{w}_0 + \sin(u \cdot \theta_w) \widetilde{w}_1}{\sin \theta_w} .
\end{align*}
\item Geodesic distance between the points is given by $\sqrt{\theta_z^2 + \theta_w^2}$.
\end{enumerate}

\subsubsection{Geodesics in the Shape Space of Framed Loops}\label{sec:geodesics_framed_loop_space}

We approximate geodesics in the shape space $\mathcal{S}_c$ as follows. For framed loops $(\gamma_0,V_0)$ and $(\gamma_1,V_1)$ of the same mod-2 linking number, let $q_0=(z_0,w_0)$ and $q_1=(z_1,w_1)$ denote lifts in the Stiefel manifold. We align the quaternionic paths according to the SVD algorithm described above. Next we approximate an optimal reparameterization by approximating the solution to \eqref{eqn:approximate_over_diff} via dynamic programming. This is combined with an optimal seed search over the pure rotations $S^1 \subset \mathrm{Diff}^+(S^1)$ to obtain an approximate optimal reparameterization in $\mathrm{Diff}^+(S^1)$. The alignment and reparameterization procedures are iterated until a stopping condition is met, then the geodesic is given explicitly in $\Gr$ by the procedure described in the previous section. Mapping the geodesic pointwise to framed curve space under $\mathrm{H}$ produces the geodesic joining the equivalence classes of $(\gamma_0,V_0)$ and $(\gamma_1,V_1)$ in $\mathcal{S}_c$.  Examples of geodesics between closed curves are provided in Section \ref{sec:applications}.

Just as in the setting of open framed curves, geodesics in the Grassmannian can leave the open submanifold $\mathrm{Gr}_2^\ast(\mathcal{V})$. This corresponds to singular framed curves along the geodesic in closed loop space. Moreover, closed loop space is disconnected, but one can construct a geodesic joining elements of different components by taking a path in the Grassmannian $\mathrm{Gr}_2(\mathcal{P}\C)$.

\section{Frame Twisting}\label{sec:frame_twisting}

\subsection{The Frame Twisting Action}

The Lie groups involved in the quotient constructions of the previous sections are by now fairly standard in the elastic shape analysis literature. On the other hand, there is another Lie group action which is genuinely unique for framed space curves. This is the action of $\mcP S^1$ by \emph{frame-twisting}. 

We represent elements of $\mathcal{P}S^1$ by $e^{i\psi}$, where $\psi$ is a real-valued function, specified up to global addition of a multiple of $2\pi$. Then a path $e^{i\psi}$ acts on a framed curve $(\gamma,V) \in \hatS_o$ via rotating each vector $V(t)$ in the normal plane to $\gamma(t)$ by the angle $\psi(t)$, according to the right hand rule. The path $e^{i\psi}$ also acts on $(z,w) \in \mathcal{P}\C^2$ by pointwise scalar complex multiplication under the identification $S^1 = \mathrm{U}(1)$ and these actions satisfy the equivariance property
\begin{equation}\label{eqn:loop_space_equivariance}
\mathrm{H}(e^{i\psi} \cdot (z,w)) = (e^{2i\psi}) \cdot \mathrm{H}(z,w).
\end{equation}
Clearly, these actions preserve pointwise norm (for quaternionic paths) and total length (for framed curves) and so descend to equivariant actions on $S_{\sqrt{2}}^\ast$ and $\hatS_0/\R^+$. It is immediately obvious that the action on $\mathcal{P}\C^2$ is by $L^2$ isometries and it follows that the action on framed loop space is isometric with respect to $g^\mathcal{S}$.

Similarly, the loop space $\mathcal{L}S^1$ acts isometrically on $\mathcal{L}\C^2 \cup \mathcal{A}\C^2$ and on the space of framed loops and this action is equivariant under the frame-Hopf map $\mathrm{H}$ in the sense of \eqref{eqn:loop_space_equivariance}.

\subsection{Principal Bundle Structures}

\subsubsection{Open Curves}

The space $\widehat{\mathcal{S}}_o$ of open framed curves is a principal bundle over the space $\mathrm{Imm}_0(I,\R^3)$ of immersed (unframed) curves based at the origin, with structure group $\mathcal{P}S^1$. Local trivializations can be built in terms of cross-sections as follows. Consider the subset $\mathcal{U} \subset \mathrm{Imm}_0(I,\R^3)$ containing curves $\gamma$ with $\gamma'(0) \neq (1,0,0)$. A cross-section is defined by the map $
\mathcal{U} \rightarrow \widehat{\mathcal{S}}_o$ taking $\gamma$ to $(\gamma,V)$ where $V$ is the unique vector field with no intrinsic twisting around $\gamma$ and such that $V(0)$ is the normalized projection of $(1,0,0)$ to the orthogonal complement of $\gamma'(0)$---framings without intrinsic twisting are known as \emph{rotation minimizing frames} \cite{wang2008computation} or \emph{Bishop frames}, after \cite{bishop1975there}. Similar charts can be used to cover all of $\mathrm{Imm}_0(I,\R^3)$.

\subsubsection{Closed Curves}

Similarly, each component of the space $\widehat{\mathcal{S}}_c$ of closed framed curves is an $\mathcal{L}S^1$-principal bundle over the space $\mathrm{Imm}_0(S^1,\R^3)$ of immersed loops based at the origin. The local trivializations are more complicated, since a rotation minimizing frame will not necessarily close. Moreover, the standard Frenet framing is not well-defined for closed curves with vanishing curvature. Nonetheless, local cross-sections can be constructed by adapting the rotation minimizing frame idea. The principal bundle structure of 
$$
\mathcal{L}S^1 \hookrightarrow \widehat{\mathcal{S}}_c \rightarrow \mathrm{Imm}_0(S^1,\R^3)
$$
is described in detail in \cite{needham2017kahler}.

\subsection{Optimal Frame Registration}

\subsubsection{Frame Registration Theorem}

The elastic shape analysis framework proposes that a geodesic between a pair of unparameterized curves can be realized as a geodesic in the total space of parameterized curves which is horizontal with respect to the action of the reparameterization group. The takeaway of the discussion in the previous section is that a geodesic between a pair of unframed curves can be obtained by finding a geodesic between framed curves which is horizontal with respect to the frame-twisting action. Finding such a horizontal geodesic once again involves an optimization over an infinite-dimensional group. Unlike the parameterization registration, the problem of optimization over $\mathcal{P}S^1$ or $\mathcal{L}S^1$ can be solved explicitly, at least for generic curves.

More precisely, given $q_0=(z_0,w_0)$ and $q_1=(z_1,w_1)$ in $S_{\sqrt{2}}^\ast$, we seek $\widehat{q}_1$ in the $\mathcal{P}S^1$--orbit of $q_1$ so that the geodesic $q_u$ joining $q_0$ and $\widehat{q}_1$ is $\mathcal{P}S^1$--horizontal. This means that the tangent $\frac{d}{du} q_u$ contains no component in the $\mathcal{P}S^1$--orbit direction for all $u$. The image of this geodesic under $\mathrm{H}$ is a path of framed curves $(\gamma_u,V_u)$ such that the variation fields along the geodesic are always orthogonal to the $\mathcal{P}S^1$-orbits. Thus $\gamma_u$ represents a geodesic in the space $\hatS_o/(\R^+ \times \mathcal{P}S^1)$ of unframed curves. 

\begin{thm}\label{thm:frame_twisting}
Let $q_0=(z_0,w_0)$ and $q_1=(z_1,w_1)$ be elements of  $\in S_{\sqrt{2}}^\ast$ which do not lie in the same $\mathcal{P}S^1$--orbit and such that $\left<(z_0(t),w_0(t)),(z_1(t),w_1(t))\right>_{\C^2} \neq 0$ for all $t \in I$. Then the spherical geodesic joining $q_0$ and 
\begin{equation}\label{eqn:optimal_framing}
\widehat{q}_1 = \frac{\left<(z_0,w_0),(z_1,w_1)\right>_{\C^2}}{\left|\left<(z_0,w_0),(z_1,w_1)\right>_{\C^2}\right|} q_1
\end{equation}
is $\mathcal{P}S^1$--horizontal.
\end{thm}

The proof of this theorem is given in Section \ref{sec:proof_frame_twisting}. 

\subsubsection{Horizontal Geodesics for Framed Curves} 

The theorem can be rephrased as follows. Let $(\gamma_0,V_0)$, $(\gamma_1,V_1) \in \widehat{\mathcal{S}}_0/\R^+$ denote framed paths with quaternionic representatives $q_0=(z_0,w_0)$, $q_1=(z_1,w_1) \in S_{\sqrt{2}}^\ast$. Then the framed curve in the $\mathcal{P}S^1$-orbit of $(\gamma_1,V_1)$ which is closest in geodesic distance to $(\gamma_0,V_0)$ is the one corresponding to $\widehat{q}_1$, as defined by \eqref{eqn:optimal_framing}. 

\subsection{Optimal Frame Registration for Closed Framed Curves}

We wish to extend the optimal framing theorem of the previous section to closed framed curves. Let 
$$
S_{\sqrt{2},l} = \left\{ q \in \mathcal{L}\mathh \mid \|q\|_{L^2}^2 = 2 \right\}
$$
denote the \emph{Hilbert sphere} in the loop space $\mathcal{L}\mathh$. Similarly, let $S_{\sqrt{2},a}$ denote the Hilbert sphere in the antiloop space. For $x=l$ or $a$, we denote by $S_{\sqrt{2},x}^\ast$ the open submanifold of  $S_{\sqrt{2},x}$ containing quaternionic curves which do not pass through the origin. It is easy to see that these manifolds are invariant under the action of the loop group $\mathcal{L}S^1$. Moreover, going through the proof of Theorem \ref{thm:frame_twisting} (see Section \ref{sec:proof_frame_twisting}), one sees that it applies to these spheres as well. We record this in the following corollary.

\begin{cor}\label{cor:optimal_framing_closed}
Let $q_0=(z_0,w_0)$ and $q_1=(z_1,w_1)$ denote elements of $S_{\sqrt{2},x}^\ast$, for $x=l$ or $a$, which do not lie in the same $\mathcal{L}S^1$--orbit and such that $\left<(z_0(t),w_0(t)),(z_1(t),w_1(t))\right>_{\C^2} \neq 0$ for all $t \in I$. Then the spherical geodesic joining $q_0$ and 
\begin{equation}\label{eqn:optimal_framing_closed}
\widehat{q}_1 = \frac{\left<(z_0,w_0),(z_1,w_1)\right>_{\C^2}}{\left|\left<(z_0,w_0),(z_1,w_1)\right>_{\C^2}\right|} q_1
\end{equation}
is $\mathcal{L}S^1$--horizontal.
\end{cor}

\begin{figure*}[t]
\begin{center}
\includegraphics[width=\textwidth]{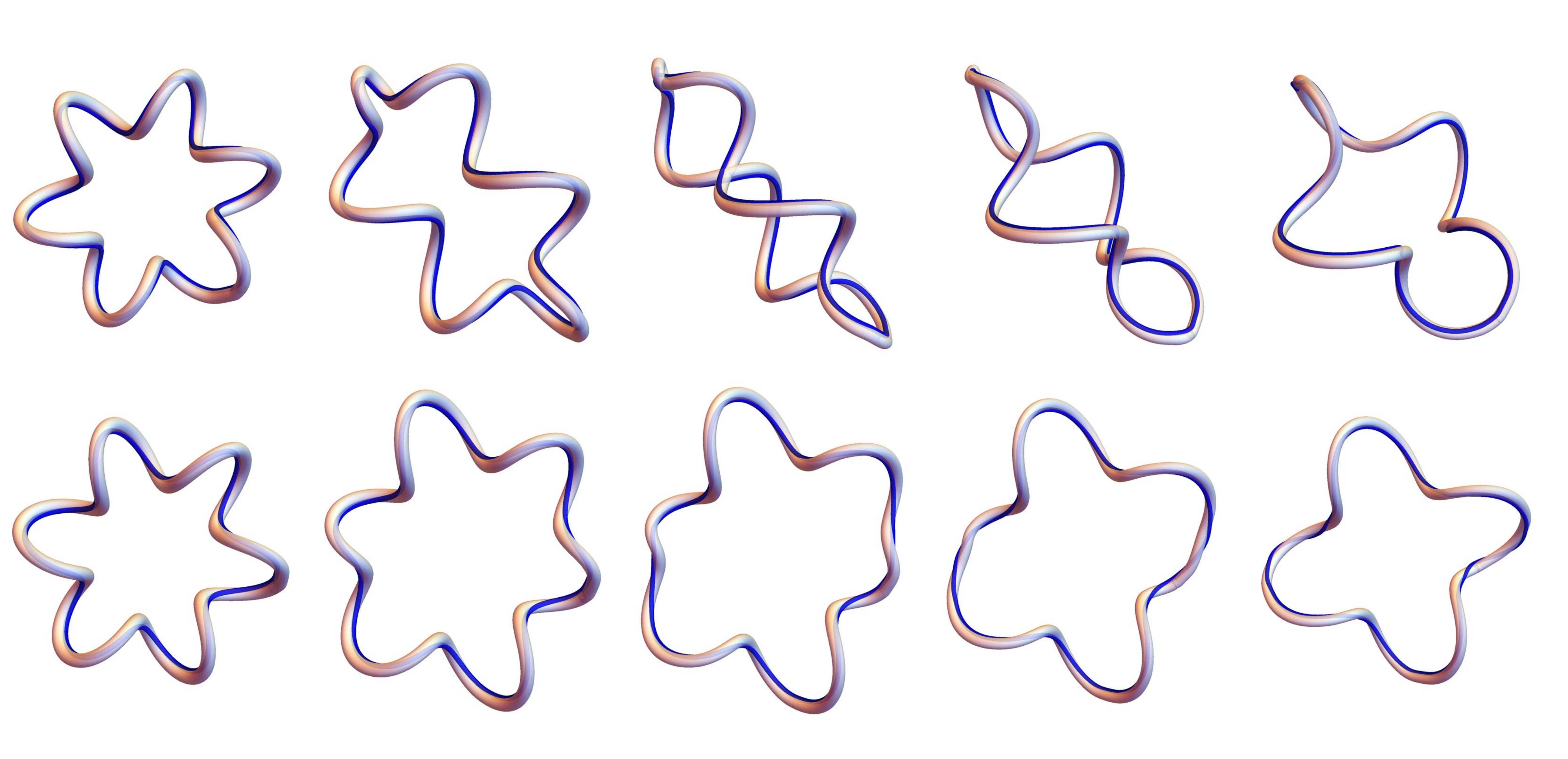}
\end{center}
\caption{Geodesic between framed spirals on a torus of revolution. The framed curves are visualized as a thickened base curve, together with a blue line along the surface indicating the twisting of the framing. In the example in the first row, the endpoints of the geodesic are endowed with their Frenet frames and the geodesic is computed via the algorithm in Section \ref{sec:geodesics_framed_loop_space}. In the example in the second row, the curve at the beginning of the geodesic has its Frenet frame, while the framing on the second curve is optimized with respect to the first curve, using the algorithm in Section \ref{sec:explicit_geodesics_for_space_curves}, and in particular  Corollary \ref{cor:optimal_framing_closed}. The  Geodesic distance in the first example is $0.29$ and geodesic distance in the second example is $0.18$---distances have been normalized so that the Grassmannian has maximum diameter equal to one. The distance in the second example represents geodesic distance in the space of unframed curves.} \label{fig:optimized_frame_geodesic}
\end{figure*}

The Stiefel manifold $\mathrm{St}_2(\mathcal{V})$ is a submanifold of the appropriate Hilbert sphere $S_{\sqrt{2},x}$ (with $x=l$ for $\mathcal{V} = \mathcal{L}\C$ and $x=a$ for $\mathcal{V} = \mathcal{A}\C$). Furthermore, the submanifold is invariant under the action of $\mathcal{L}S^1$. We therefore have the following interpretation of the corollary.  Let  $(\gamma_0,V_0)$, $(\gamma_1,V_1)$ denote framed loops of the same mod-2 linking number, with quaternionic representatives $q_0=(z_0,w_0)$, $q_1=(z_1,w_1) \in \mathrm{St}_2(\mathcal{V})$. The quaternionic curve in the $\mathcal{L}S^1$-orbit of $q_1$ which is closest to $q_0$ with respect to geodesic distance in $S_{\sqrt{2},x}$ is $\widehat{q}_1$, as defined by \eqref{eqn:optimal_framing_closed}. Since $\mathrm{St}_2(\mathcal{V})$ is an $\mathcal{L}S^1$-invariant submanifold of $S_{\sqrt{2},x}$ with the induced Riemannian metric, $\widehat{q}_1$ also minimizes distance to $q_0$ over the $\mathcal{L}S^1$-orbit in the Stiefel manifold. It follows that the framed loop in the $\mathcal{L}S^1$-orbit of $(\gamma_1,V_1)$ which is closest in geodesic distance to $(\gamma_0,V_0)$ is the one corresponding to $\widehat{q}_1$. See Figure \ref{fig:optimized_frame_geodesic} for an example of the effect of optimizing the frame.

\subsection{Geodesics for Unframed Curves}

We can incorporate the frame optimization idea described above into our algorithms to produce new algorithms for approximating geodesics between unframed curves. The case of closed curves will be described in detail here, but geodesics for open curves can be treated similarly. 

\subsubsection{Geodesics for Closed Space Curves}\label{sec:explicit_geodesics_for_space_curves}

Let $\gamma_0$ and $\gamma_1$ denote parameterized closed space curves, preprocessed so that $\mathrm{length}(\gamma_j) = 2$. Choose a framing $V_j$ for $\gamma_j$ so that the resulting framed curves have the same mod-2 linking number (e.g., start with Frenet framings or minimal twist framings of \cite{wang2008computation}, then adjust as necessary to meet the linking number condition). Let $q_j=(z_j,w_j)$ denote the Stiefel manifold representation of $(\gamma_j,V_j)$. Then
\begin{enumerate}
\item the SVD alignment of Section \ref{sec:grassmannian_geodesics} is applied to produce $\mathrm{U}(2)$-aligned curves $q_0^\ast$ and $q_1^\ast$,
\item a distance-minimizing reparameterization $\widehat{q}_1$ of $q_1^\ast$ is approximated via dynamic programming, and 
\item the optimal point $\widetilde{q}_1$ in the $\mathcal{L}S^1$--orbit of $\widehat{q}_1$ is computed explicitly via Corollary \ref{cor:optimal_framing_closed}. 
\end{enumerate}
This procedure is iterated until a stopping condition is met to produce optimally parameterized and framed curves. The explicit geodesic joining the resulting curves in $\Gr$ is computed via the procedure described in Section \ref{sec:grassmannian_geodesics}, and applying the frame-Hopf map at each point along the geodesic produces a geodesic $\left(\widetilde{\gamma}_u,\widetilde{V}_u\right)$ between framed curves. The homotopy of the base curve $\widetilde{\gamma}_u$ represents a geodesic in the space of closed curves between the $(\mathrm{SO}(3) \times \mathrm{Diff}^+(S^1))$--orbits of the initial curves $\gamma_0$ and $\gamma_1$.

\section{Numerical Examples}\label{sec:applications}

\subsection{Geodesic Examples}

\begin{figure*}
\begin{center}
\includegraphics[width=\textwidth]{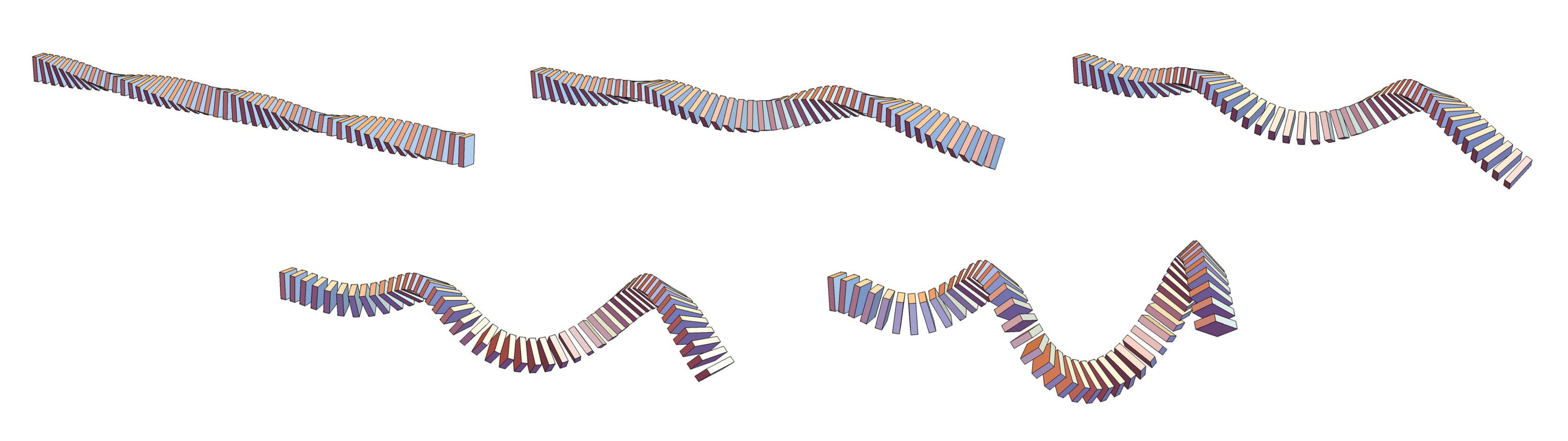}
\end{center}
\caption{A geodesic between framed curves, illustrating an application to trajectory blending. Each framed curve is visualized as a spatiotemporal trajectory of a rigid object. The starting point (top left) is a ``barrel roll" trajectory, the finishing point (bottom right) is a more complicated trajectory, and intermediate points along the geodesic are blended trajectories. Geodesic representing trajectory blending. Geodesic distance is 0.16, where the space of open curves (i.e., the Hilbert sphere) has been normalized to have diameter one.}\label{fig:trajectory_blending}
\end{figure*}

Figure \ref{fig:trajectory_blending} shows a basic example of a geodesic between open framed curves, visualized in the context of  \emph{trajectory blending}. A common animation technique, known as \emph{motion blending}, is to produce a new animation via interpolation between existing ones---see \cite{kovar2003flexible} for a general introduction and \cite{bauer2017landmark,celledoni2016shape} for applications of the elastic shape analysis framework to the topic. Adapted framed curves give a natural way to describe trajectories of rigid bodies such as aircraft, and our results provide a method of blending such trajectories.

\begin{figure*}
\begin{center}
\includegraphics[width=\textwidth]{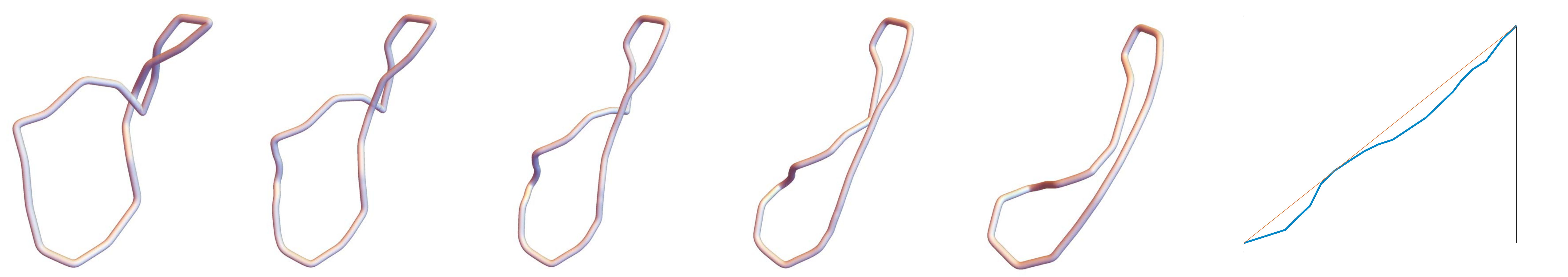}
\end{center}
\caption{A geodesic in the space of (unframed) closed curves between circular DNA molecules. The geodesic is calculated by initializing the algorithm with Rotation Minimizing Frames, then optimizing the framing for the ending curve. The geodesic distance between the molecules is $0.21$, where the space has been rescaled to have diameter 1. The curves were initialized with arclength parameterization, then the optimal reparameterization of the ending curve was approximated via dynamic programming. The plot on the right shows the optimal reparameterization. It illustrates the matching of common features between the curves at the endpoints of the geodesic.} \label{fig:DNA_minicircles_geodesic}
\end{figure*}

Figure \ref{fig:DNA_minicircles_geodesic} shows a geodesic in the space of (unframed) closed curves, using the algorithm of Section \ref{sec:explicit_geodesics_for_space_curves}. The endpoints of the geodesic are curves representing closed DNA molecules. See Section \ref{sec:DNA_minicircles} for a description of the dataset these curves were taken from.

\begin{figure*}
\begin{center}
\includegraphics[width=0.9\textwidth]{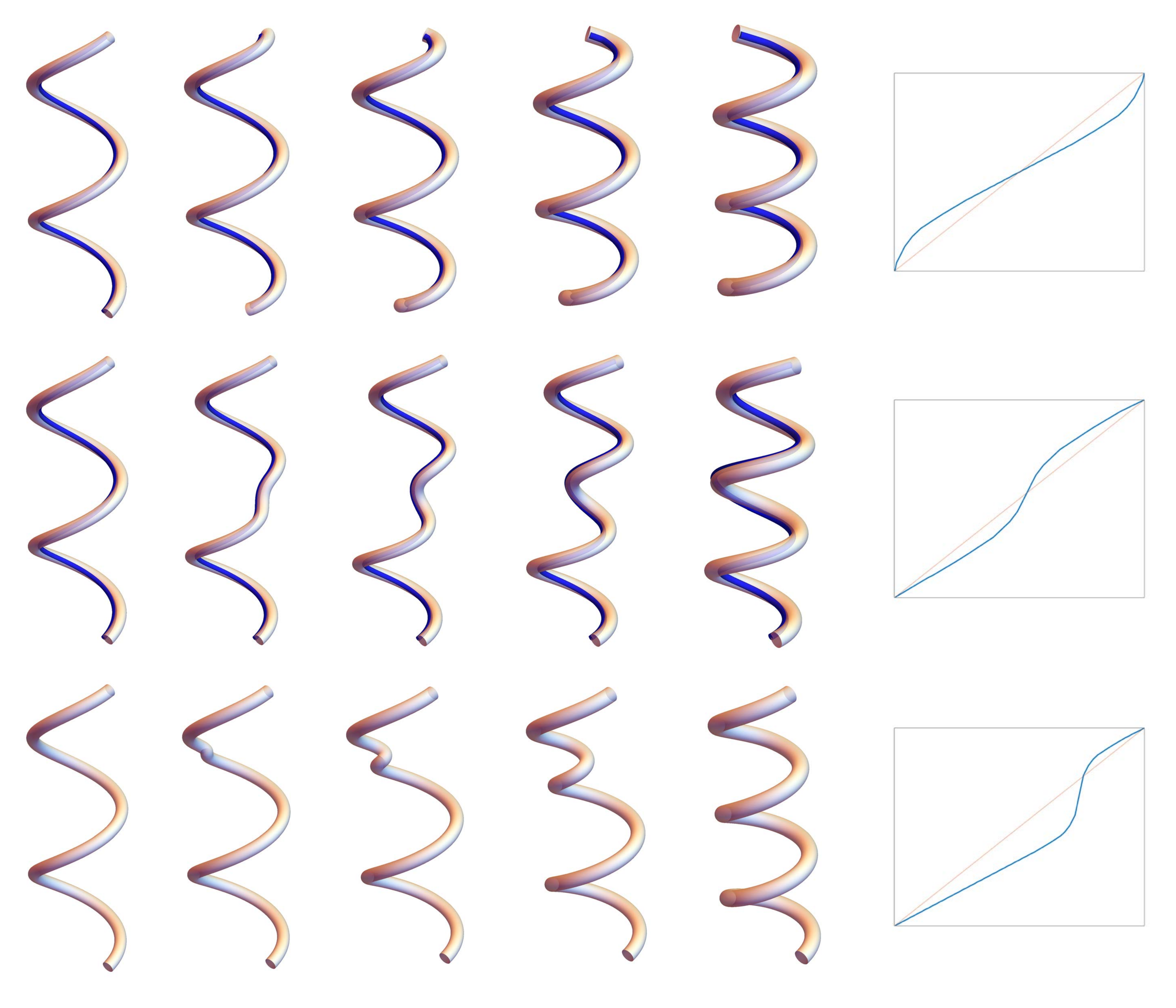}
\end{center}
\caption{Geodesics between a helix of 2 turns and a helix of 3 turns. Each row shows curves evenly spaced along the geodesic so that the curve in the middle is Fr\'{e}chet mean of the endpoints. The last figure in each row shows the optimal reparameterization (approximated via dynamic programming) of the 3 turn helix, which was initialized with arclength parameterization. The geodesics in the first two rows are in framed curve space and framed curves are visualized as thickened base curves with a blue line on their surfaces indicating twisting of the framing. The first row shows a geodesic where the helices on the endpoints are endowed with their Frenet framings. The second row shows a geodesic where the helix on the left is given its Frenet framing, and the helix on the right has its optimal framing given by Theorem \ref{thm:frame_twisting}, and the figure therefore represents a geodesic in (unframed) curve space. For comparison, the geodesic in the third row is computed according to the standard SRVT algorithm.}\label{fig:helix_geodesics}
\end{figure*}

\begin{figure*}
\begin{center}
\includegraphics[width=\textwidth]{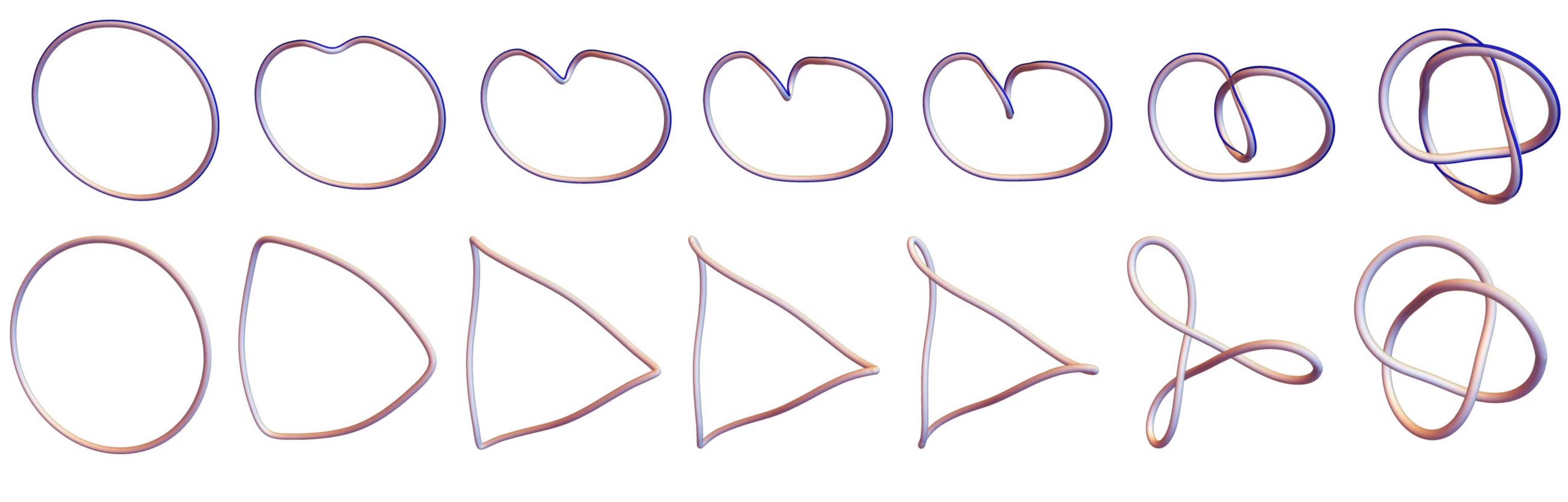}
\end{center}
\caption{Geodesics between a circle and a trefoil knot in the space of closed curves. In the geodesic in the top row, the circle is endowed with its Frenet frame and the trefoil is endowed with the corresponding optimal frame. For comparison, the geodesic between the curves with respect to the SRVT algorithm is shown in the bottom row. Note that the SRVT geodesic passes through a singular curve, while the framed curve geodesic does not.}\label{fig:circle_torus_geodesics}
\end{figure*}

More examples of geodesics between open and closed framed curves are provided in Figures \ref{fig:helix_geodesics} and \ref{fig:circle_torus_geodesics}. For comparison, we also compute geodesics via the SRVT algorithm \footnote{In particular, we used the Matlab implementation available on the FSU Statistical Shape Analysis and Modeling Group website http://ssamg.stat.fsu.edu} of \cite{joshi2007novel,srivastava2011shape}.

\subsection{DNA Minicircles}\label{sec:DNA_minicircles}

In \cite{irobalieva2015structural}, the authors study short (336 bp), circular DNA molecules, referred to as \emph{DNA minicircles}. The shape data from their experiment consists of parameterized space curves representing the centerlines of DNA minicircle double-helices, obtained via electron cryo-tomography. There are nine subcollections of such space curves, corresponding to different degrees of supercoiling in the molecules. For example, one of the subcollections contains samples of \emph{relaxed} DNA molecules. In its relaxed state, the two strands of a DNA minicircle wrap around one another 32 times to form the familiar double-helix structure. The remaining subcollections contain DNA molecules with different numbers of wrappings of the two strands; we refer to this distinguishing feature as a \emph{link defecit}. The subcollections provided to us contain space curves corresponding to DNA minicircles with link defecits $\Delta \mathrm{Lk} = - 6$, $-4$,  $-3$, $-2$, $-1$, $0$, $1$, $2$ and $3$, respectively. In \cite{irobalieva2015structural}, it was observed by human classification of the space curves into various shape motifs that link defecit plays a role in the shapes of the centerline curves. Roughly, larger link defecits correspond to molecule centerlines which are more compact and twisted. 

\begin{figure*}
\begin{center}
\includegraphics[width=\textwidth]{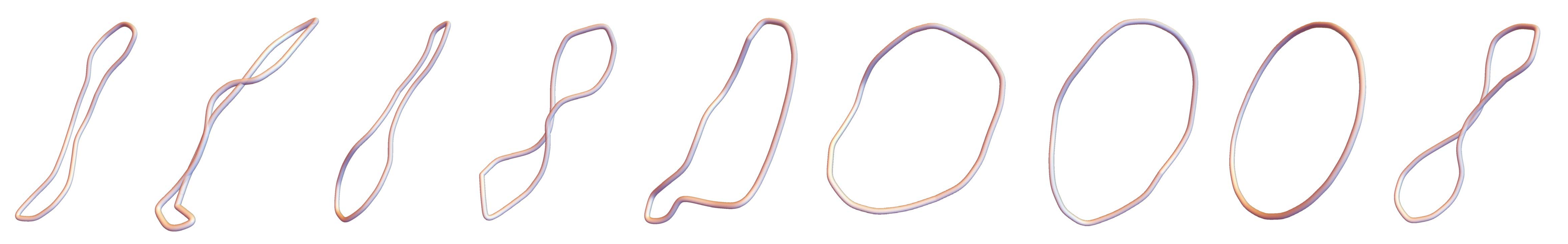}
\end{center}
\caption{Centers of the largest medoids for each collection of DNA minicircles in the $k$-medoid clustering experiment. From left to right, the centers correspond to DNA minicircles with $\Delta \mathrm{Lk} = -6, -4, -3, -2, -1, 0, 1, 2, 3$. The percentage of samples contained in the largest medoid for the $\Delta \mathrm{Lk} = -6$ collection is $25\%$. The respective percentages for the remaining collections are  $38.9\%, 39.1\%, 35.0\%, 35.0\%, 46.7\%, 40.4\%, 37.9\%, 29.4\%$. This figure illustrates the qualitative description of shape variation with respect to $\Delta \mathrm{Lk}$ given in \cite{irobalieva2015structural} .}\label{fig:minicircle_medoids}
\end{figure*}

To get a better picture of the shape variation accross link defecit categories, we performed the following clustering experiment. We computed the distance matrix within each $\Delta \mathrm{Lk}$ class, then ran a $k$-medoid clustering algorithm with $k=8$ (the authors of \cite{irobalieva2015structural} classify DNA minicircles into 8 shape categories). The centers of the largest medoid clusters for each link defecit category are shown in Figure \ref{fig:minicircle_medoids}. This picture agrees qualitatively with the conclusions of \cite{irobalieva2015structural}.

\subsection{Flag Means}

In the context of shape classficiation, a useful calculation is the Karcher mean of a collection of shapes. The standard procedure in the elastic shape analysis framework for computing Karcher means is an iterative algorithm utilizing the logarithm map to do calculations in linear tangent spaces. There is an extra numerical step in this approach, because the logarithm and exponential maps for the submanifold of closed curves do not have an explicit form.

In the setting presented here, computations can be done explicitly on the submanifold of closed curves, as spaces of closed (planar or framed space) curves are identified with Grassmannians. There is a recently proposed alternative notion of average for a collection of samples of a Grassmann manifold called the \emph{flag mean} \cite{draper2014flag}. Let $\{[q_j]=[z_j,w_j]\}$, $j=1,\ldots,N$, be a collection of points in one of the Grassmannians $\mathrm{Gr}_2(\mathcal{L}\R)$, $\mathrm{Gr}_2(\mathcal{L}\C)$ or $\mathrm{Gr}_2(\mathcal{A}\C)$. The flag mean $[q]=[z,w]$ of the samples is given by first finding
$$
\mathrm{argmin}_{[z]} \sum_{j=1}^N d_{pF}([z],[q])^2,
$$  
where $[z]$ is the one-dimensional subspace spanned by $z$ and $d_{fP}$ is projection-Frobenius distance. Next $[w]$ is obtained as
$$
\mathrm{argmin} \sum_{j=1}^N d_{pF}([w],[q])^2,
$$
where the minimum is taken over $w$ such that $[w]$ is orthogonal to $[z]$.  While the flag mean is not, in general, equal to the Karcher mean, it does give a good approximation for well-clustered data. Moreover, there is a simple algorithm for computing the flag mean in finite-dimensional Grassmannians \cite{draper2014flag} which extends without issue to our infinite-dimensional setting. 

We choose to utilize the flag mean in our algorithm for computing curve averages, since it allows us to work explicitly and directly on the submanifold of closed curves as often as possible. We will describe the details of our curve averaging algorithm in the case of closed space curves. Similar algorithms work for closed plane curves and for closed framed curves in $\R^3$.

Let $\{\gamma_j\}$, $j=1,\ldots,N$ be a collection of closed curves, preprocessed to have length $2$. Choose  an arbitrary framing $V_j$ for each $\gamma_j$ so that all framed curves have the same mod-2 linking number and let $q_j = (z_j,w_j)$ be a Stiefel manifold representative of $(\gamma_j,V_j)$. We initialize with $q^\ast = q_1$. For each $j$, align $q_j$ with $q^\ast$ via SVD,  find an optimal reparameterization $\widehat{q}_j$ with respect to $q^\ast$ (approximated via dynamic programming), then find an optimal element $\widetilde{q}_j$ of the $\mathcal{L}S^1$--orbit of $\widehat{q}_j$ with respect to $q_1$ (this is done explicitly via \ref{cor:optimal_framing_closed}). Iterate this subloop for each $j$ until a stopping condition is met to obtain a collection of $\widetilde{q}_j$ which have been aligned over rotations, reparameterizations and framings to $q^\ast$. Next set $q^\ast$ to be the flag mean of the samples $\widetilde{q}_j$ and iterate the whole procedure until a stopping condition is met. The final flag mean is then mapped to the average closed curve $\gamma^\ast$ via the first coordinate of $\mathrm{H}$.

\begin{figure*}
\begin{center}
\begin{tabular}{|c|c|}
\hline
Samples & Flag Mean \\
\hline 
\includegraphics[scale=0.15]{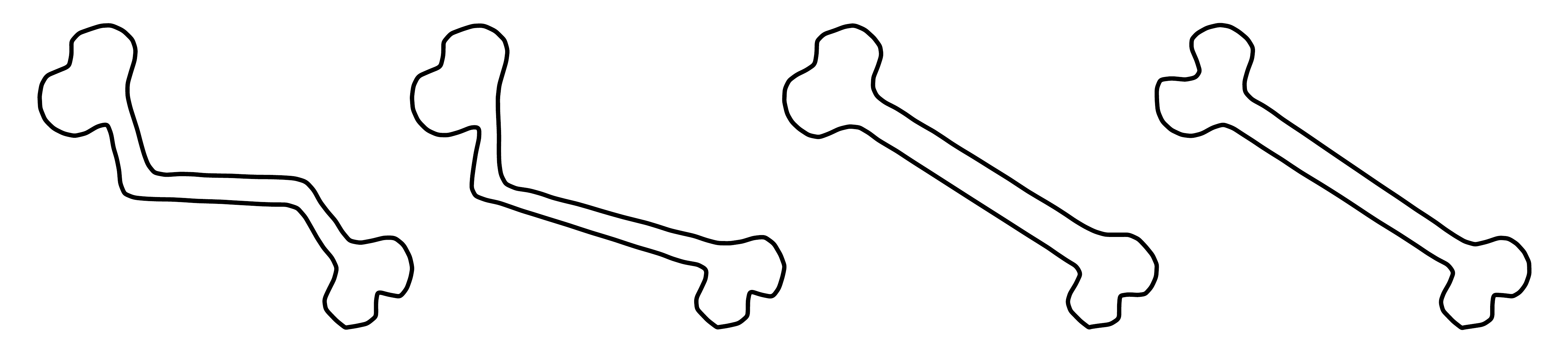} & \includegraphics[scale=0.15]{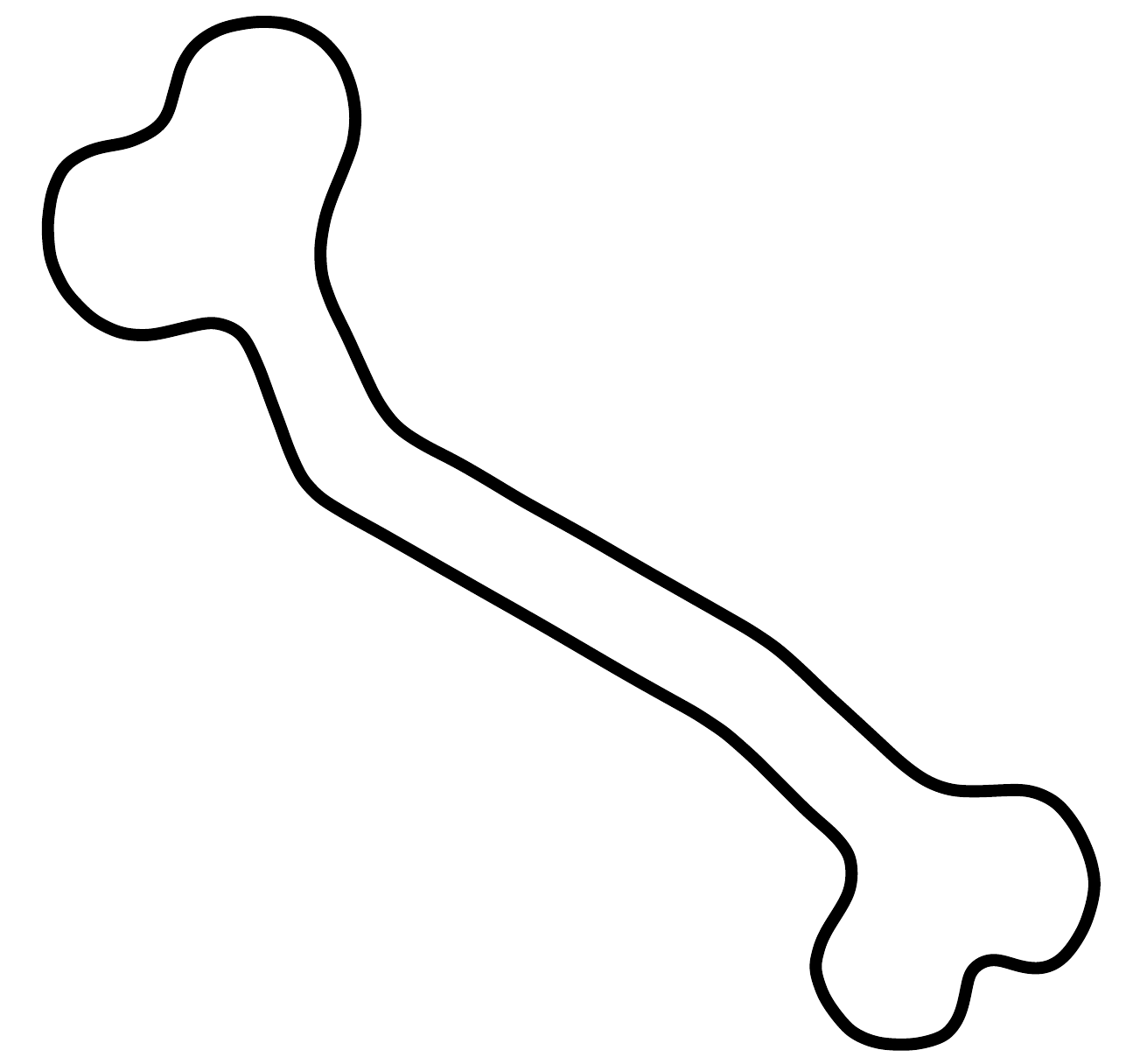} \\
\hline 
\includegraphics[scale=0.15]{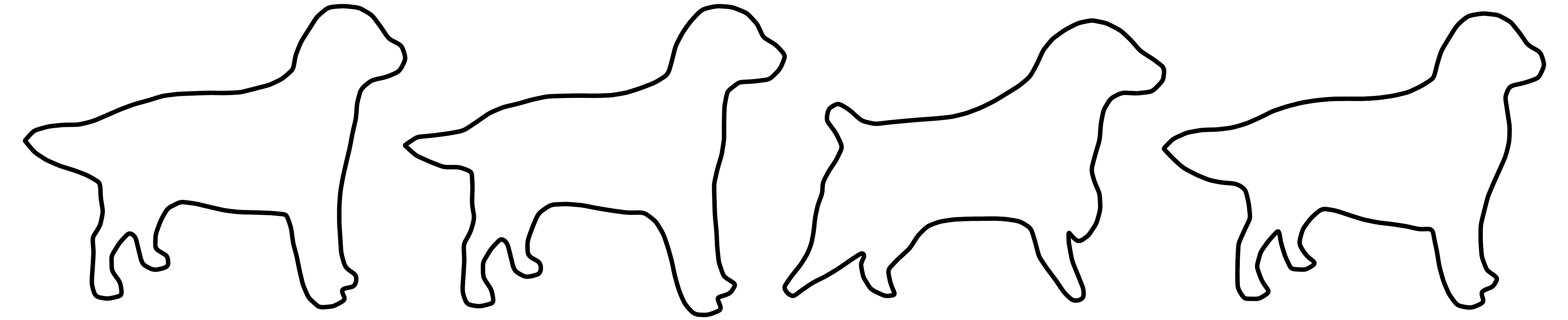} & \includegraphics[scale=0.15]{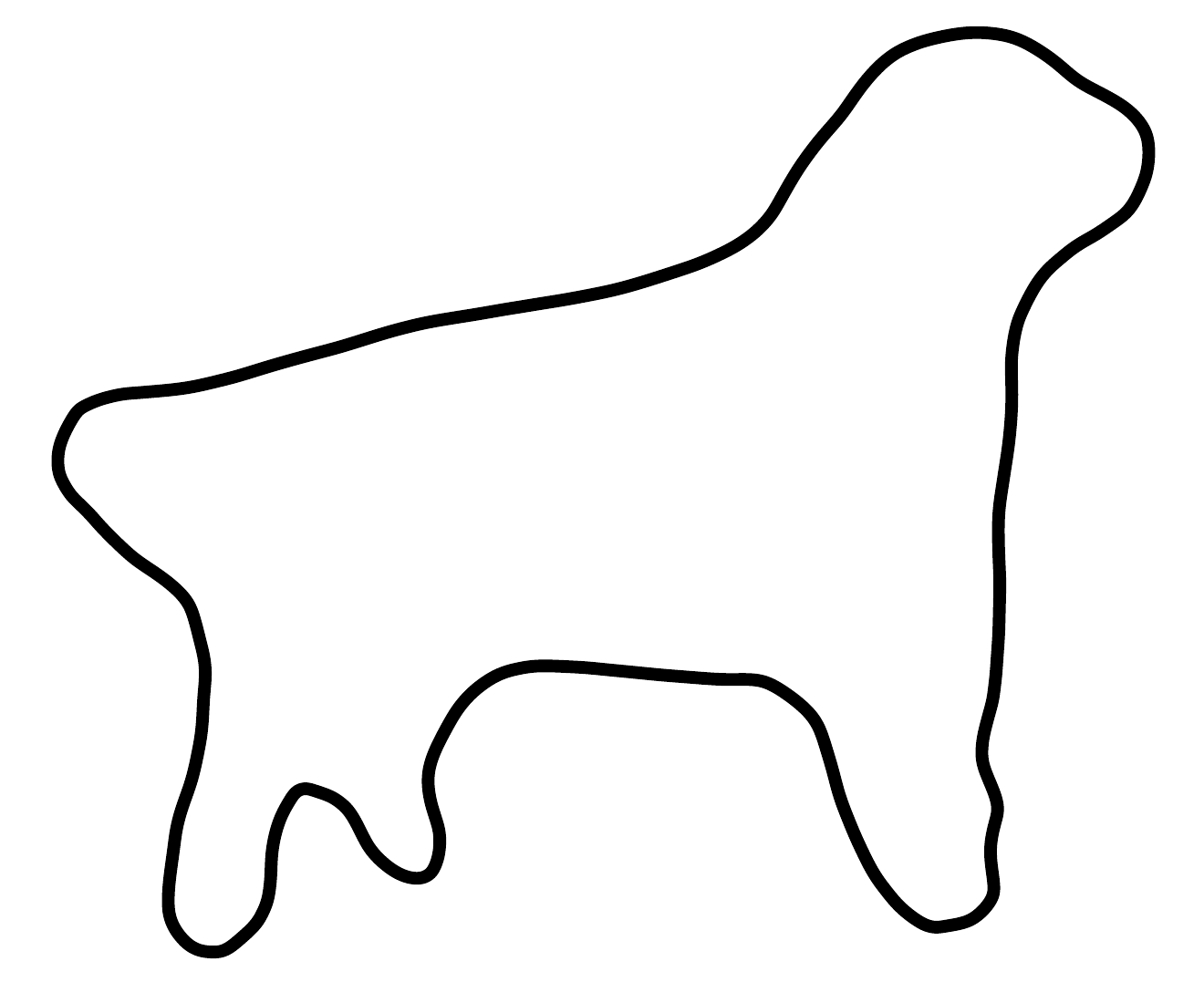} \\
\hline
\includegraphics[scale=0.15]{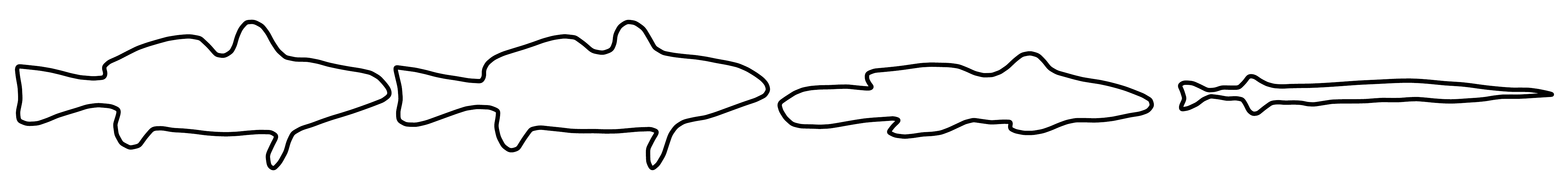} & \includegraphics[scale=0.15]{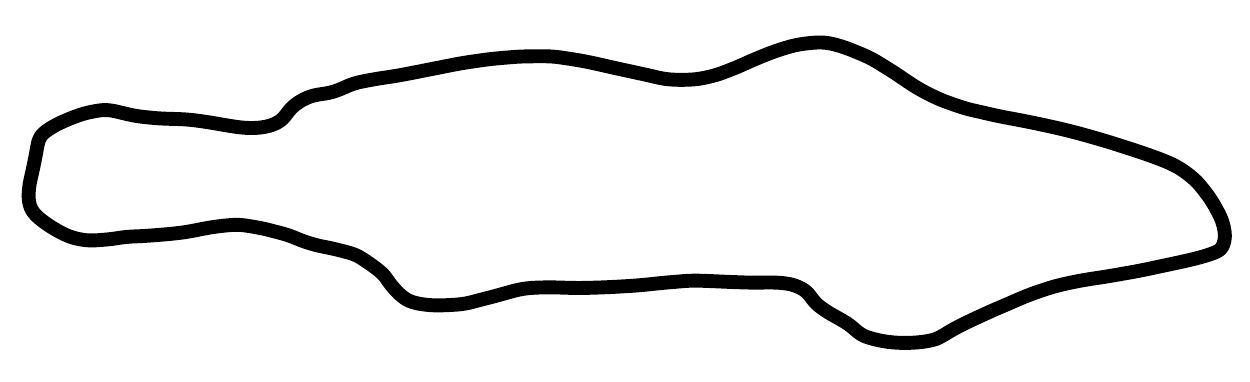} \\
\hline
\end{tabular}
\end{center}
\caption{Averages for collections of planar curves computed according to the Flag Mean algorithm. Curve samples are taken from the popular MPEG7 computer vision shape database.} \label{fig:flag_means_planar}
\end{figure*}

\begin{figure*}
\begin{center}
\begin{tabular}{|c|c|}
\hline
Samples & Flag Mean \\
\hline 
\includegraphics[scale=0.35]{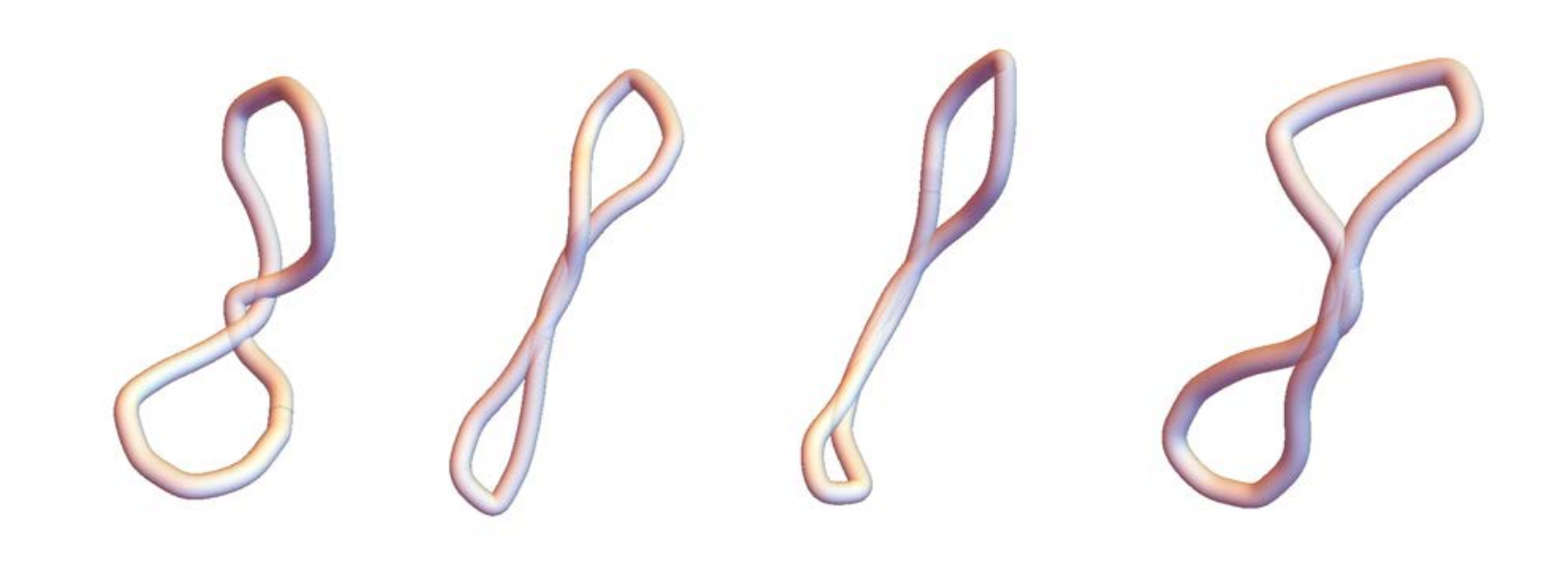} & \includegraphics[scale=0.35]{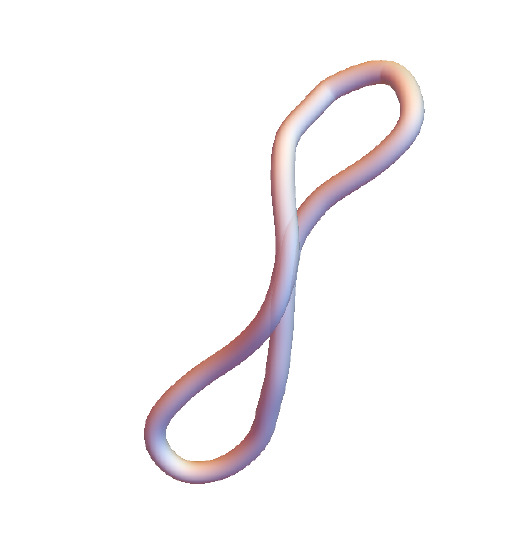} \\
\hline 
\includegraphics[scale=0.35]{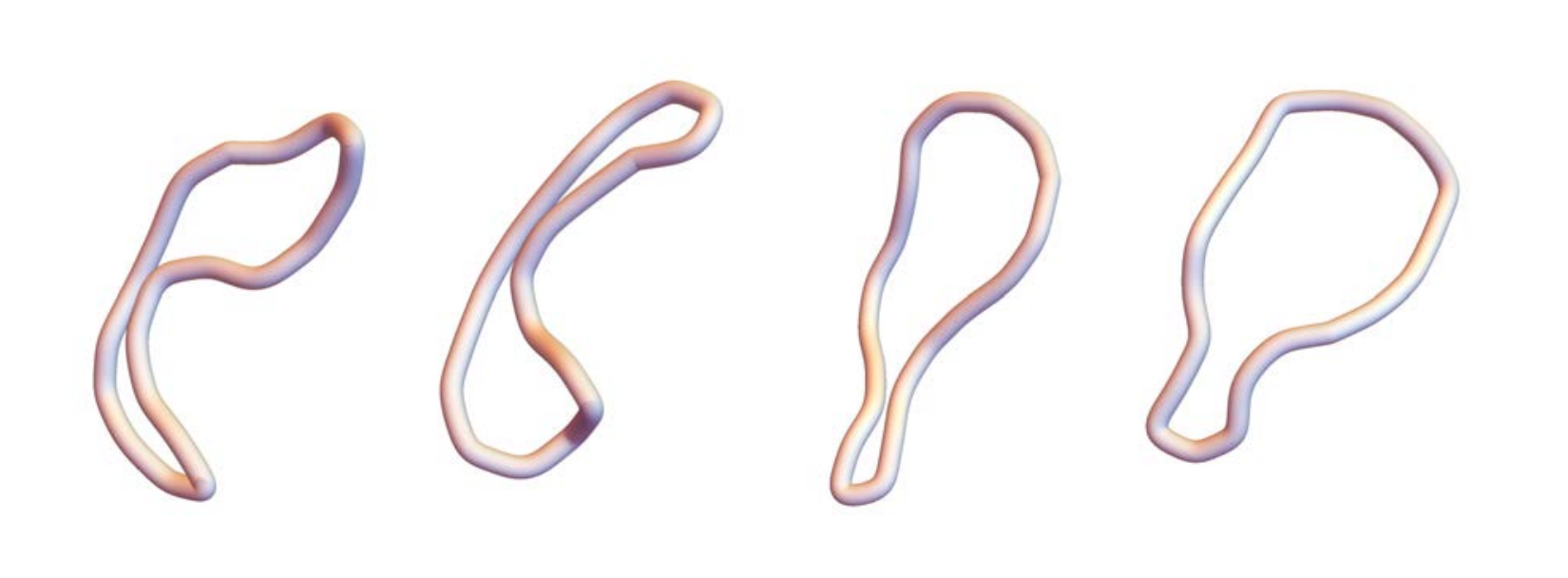} & \includegraphics[scale=0.35]{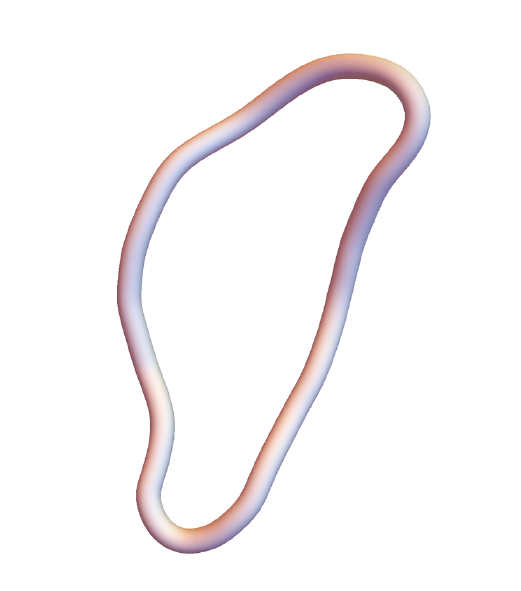} \\
\hline 
\end{tabular}
\end{center}
\caption{Averages of DNA minicircle samples computed via the Flag Mean algorithm.}\label{fig:flag_means_DNA}
\end{figure*}

Examples of flag means for plane curves and DNA minicircles are shown in Figures \ref{fig:flag_means_planar} and \ref{fig:flag_means_DNA}, respectively.

\section{Discussion}

There are several interesting directions in which to pursue further research. In \cite{kurtek2018simplifying}, the square root transform of Younes et.\ al.\ is extended to give simplifications of all elastic metrics on open planar curves. The extended transforms once again leverage the ability to perform complex arithmetic pointwise. It seems likely that a similar generalization could be performed in the framed curve setting using quaternionic arithmetic to provide simplifications of all framed curve elastic metrics. Another avenue for future research would be to develop algorithms for statistics on the space of closed (planar or space) curves by utilizing the Grassmannian correspondence, as we have begun to do here by utilizing the flag mean. The necessary exponential maps on the relevant infinite-dimensional Stiefel manifolds have already been worked out in \cite{harms2012geodesics,needham2017kahler}. Finally, precise results on the existence of optimal reparameterizations in the SRVT framework were provided in the recent paper \cite{bruveris2016optimal}. The results there rely on the existence of explicit geodesics, and therefore only apply to open curves. It would be interesting to try to extend these results to the metrics used in this paper, since they are so far the only metrics where explicit geodesics can be described for closed curves.

\section*{Acknowledgements}

I would like to thank Muyuan Chen, Steven Ludtke and Lynn Zechiedrich for graciously providing me with the very interesting DNA minicircles data. Next I would like to thank Michael Tychonievich for his help in developing a GUI for the framed curves matching program used to produce the numerical experiments. Finally, many thanks are due to various colleagues with whom I have had conversations about elastic shape analysis and framed curves over the years, including Jason Cantarella, Sebastian Kurtek, Erik Schreyer and Clayton Shonkwiler.

\section{Appendix}

\subsection{Proof of Theorem \ref{thm:pullback_metric}}\label{sec:proof_pullback}

Let $g^{\mathrm{SO}(3)}$ denote the standard bi-invariant metric on $\mathrm{SO}(3)$ induced by the Euclidean metric on $\mathfrak{so}(3) \approx \R^3$, let $g^{\R^+}$  denote the bi-invariant metric on $\R^+$ induced by 
$$
(r_1,r_2) \mapsto r_1 r_2
$$
on $T_1 \R^+ = \R$ and let $g^{\mathrm{SO}(3)} \otimes g^{\R^+}$ denote the product metric on $\mathrm{SO}(3) \times \R^+$. A natural $L^2$-type metric on $\mathcal{P}(\mathrm{SO}(3) \times \R^+)$ is given by
$$
g_{(A,r)} (\cdot,\cdot) = \frac{1}{4} \int_0^2 \left(g^{\mathrm{SO}(3)} \otimes g^{\R^+}\right)_{(A(t),r(t))} (\cdot,\cdot) \mathrm{d}s.
$$
where $(A,r) \in \mathcal{P}(\mathrm{SO}(3) \times \R^+)$, the arguments of $g_{(A,r)}$ are elements of 
$$
T_{(A,r)} \mathcal{P}(\mathrm{SO}(3) \times \R^+)\approx \mathcal{P}\mathfrak{so}(3)\times \mathcal{P}\R
$$
and $\mathrm{d}s =  r(t)\mathrm{d}t$. It is straightforward to show that the pullback of $g$ to $\widehat{\mathcal{S}}_o$ via \eqref{eqn:map_to_SO3} is exactly the metric $g^\mathcal{S}$.

Now let $h$ denote the classical Hopf map \eqref{eqn:homomorphism} and let $\mathrm{sq}$ denote the squaring map $r \mapsto r^2$ for $r \in \R^+$. It is a classical fact that $h$ satisfies $h^\ast g^{\mathrm{SO}(3)} = 4 g^{\mathrm{SU}(2)}$, where $g^{\mathrm{SU}(2)}$ is the standard metric on $\mathrm{SU}(2)$, which is isometric to the round metric $g^{S^3}$ on $S^3 \approx \mathrm{SU}(2)$. It follows that 
\begin{align*}
(h \times \mathrm{sq})^\ast \left(g^{\mathrm{SO}(3)} \otimes g^{\R^+}\right) &= 4 g^{\mathrm{SU}(2)} \otimes g^{\R^+} \\
&= 4 g^{S^3} \otimes g^{\R^+},
\end{align*}
where $g^{S^3} \otimes g^{\R^+}$ is the product metric on $S^3 \times \R^+$. Let $f:\mathh \setminus \{\vec{0}\} \rightarrow S^3 \times \R^+$ denote the polar coordinate map $q \mapsto (q/\|q\|_{\mathh},\|q\|_{\mathh})$. An elementary computation shows 
$$
f^\ast \left(g^{S^3} \otimes g^{\R^+}\right)_{q} = \mathrm{Re}\left<\cdot,\cdot\right>_{\mathh}/\|q\|_{\mathh}^2.
$$
Note that the map $\mathrm{H}$ is obtained by applying $(h \times \mathrm{sq}) \circ f$ pointwise and then composing the result  with the inverse of \eqref{eqn:map_to_SO3}.  The proof is then concluded by the following calculation, in which $q \in \mathcal{P}\mathh^\ast$ satisfies $\mathrm{H}(q)=(\gamma,V)$ and $(\gamma,V) \mapsto (A,r)$ under \eqref{eqn:map_to_SO3}:
\begin{align*}
g^{L^2}_q  &= \int_I \mathrm{Re}\left<\cdot,\cdot\right>_{\mathh} \; \mathrm{d}t \\
&= \int_I f^\ast \left(g^{S^3} \otimes g^{\R^+}\right)_{q} \|q\|_{\mathh}^2 \; \mathrm{d}t \\
&= \int_I f^\ast (h \times \mathrm{sq})^\ast \left(g^{\mathrm{SO}(3)} \otimes g^{\R^+}\right)_{(A,r)} \; r(t) \mathrm{d}t \\
&= \frac{1}{4} \int_I \left((h \times \mathrm{sq}) \circ f\right)^\ast \left(g^{\mathrm{SO}(3)} \otimes g^{\R^+}\right)_{(A,r)} \mathrm{d}s \\
&= \mathrm{H}^\ast g^\mcS_{(\gamma,V)}
\end{align*}



\subsection{Proof of Proposition \ref{prop:rotation_optimization}}\label{sec:proof_rotations}

Since $\mathrm{SU}(2)$ acts by $L^2$ isometries, we seek the minimizer $\widehat{A}$ of $\arccos \left<q_0,q_1 \cdot A\right>_{L^2}$, which is equivalent to finding the maximizer of $\left<q_0,q_1 \cdot A\right>_{L^2}$. The latter quantity is equal to 
\begin{align*}
\mathrm{Re} \, \int_I q_0 \cdot \overline{q_1 \cdot A} \; \mathrm{d}t & = \mathrm{Re} \, \int_I q_0 \cdot \overline{A} \cdot \overline{q_1} \; \mathrm{d}t \\
&= \mathrm{Re} \, \overline{A} \cdot \int_I \overline{q_1} \cdot q_0 \; \mathrm{d}t \\
&= \left<\overline{A},\int_I \overline{q_1} \cdot q_0 \; \mathrm{d}t\right>_\mathh,
\end{align*}
where the second equality follows by cyclic permutation-invariance of the real part of quaternionic arithmetic. The quantity is therefore maximized by $\widehat{A} \in S^3 \approx \mathrm{SU}(2)$ with conjugate in the same direction as $\int_I \overline{q_1} \cdot q_0 \; \mathrm{d}t$, and this completes the proof.

\subsection{Proof of Theorem \ref{thm:frame_twisting}}\label{sec:proof_frame_twisting}

Let $q \in S_{\sqrt{2}}$. The \emph{horizonal tangent space} to $q$ is the subset of tangent vectors in 
$$
T_q S_{\sqrt{2}} = \{p \in \mathcal{P}\mathh \mid \left<q,p\right>_{L^2} = 0\}
$$
which are $L^2$--orthogonal to the $\mathcal{P}S^1$--orbit directions at $q$. These orbit directions are of the form $i \xi \cdot q$, where $\xi:\R \rightarrow \R$ is a smooth function. A tangent vector $p$ is therefore horizontal if and only if $\left<p,i\xi q\right>_{L^2} = 0$ for all $\xi$. Switching to complex coordinates $q=(z,w)$ and $p=(u,v)$, this condition becomes
\begin{align*}
0 &= \int_I \mathrm{Re} \left<(u,v), i \xi (z,w) \right>_{\C^2} \; \mathrm{d}t \\ 
&= \int_I -\xi \mathrm{Im}\left<(u,v),(z,w)\right>_{\C^2} \; \mathrm{d}t
\end{align*}
for all smooth $\xi$. By the standard argument from the calculus of variations, we conclude that $p=(u,v)$ is horizontal if and only if $\mathrm{Im} \left<(u,v),(z,w)\right>_{\C^2}$ is identically zero. 

Consider elements $q_0=(z_0,w_0)$ and $q_1=(z_1,w_1)$ of $S_{\sqrt{2}}^\ast$ which do not lie in the same $\mathcal{P}S^1$--orbit and with 
$$
\left<(z_0(t),w_0(t)),(z_1(t),w_1(t))\right>_{\C^2} \neq 0
$$
for all $t$. We seek $\widehat{q}_1 = e^{i\psi} \cdot q_1$ in the $\mathcal{P}S^1$--orbit of $q_1$ such that the geodesic $q_u$ joining $q_0$ and $\widehat{q}_1$ in $S_{{\sqrt{2}}}$ is horizontal for all $u$. Since $\mathcal{P}S^1$ acts by isometries, if the geodesic starts horizontal then it will stay horizontal---that is, it suffices to find $\widehat{q}_1$ so that $\left.\frac{d}{du}\right|_{u=0} q_u$ is $\mathcal{P}S^1$--horizontal at $q_0$. 

The geodesic joining $q_0$ and $\widehat{q}_1$ is given by \eqref{eqn:spherical_geodesic}. The derivative at $u=0$ of this geodesic is given by 
$$
-\frac{\theta \cos \theta}{\sin \theta} q_0 + \frac{\theta}{\sin \theta} \widehat{q}_1.
$$
Writing $q_0=(z_0,w_0)$, $q_1=(z_1,w_1)$ and recalling that $\widehat{q}_1=e^{i\psi} \cdot q_1$ for some $\psi:\R \rightarrow \R$, the desired horizontality condition reduces to
$$
\mathrm{Im} \, e^{i \psi} \cdot \left<(z_1,w_1),(z_0,w_0)\right>_{\C^2} = 0,
$$
and this condition is achieved by taking 
$$
e^{i\psi} = \frac{\left<(z_0,w_0),(z_1,w_1)\right>_{\C^2}}{\left|\left<(z_0,w_0),(z_1,w_1)\right>_{\C^2}\right|}.
$$


\bibliographystyle{spmpsci}

\bibliography{needham_bibliography}

\end{document}